%% file: main.tex
\documentclass[12pt]{article}
\usepackage{times}
\usepackage[a4paper, total={6in, 8in}]{geometry}
\usepackage[authoryear, round]{natbib}
\usepackage{stmaryrd} 
\SetSymbolFont{stmry}{bold}{U}{stmry}{m}{n} 
\usepackage{bm} 
\usepackage{dsfont} 
\usepackage{xcolor} 
\usepackage{adjustbox}
\usepackage{amsmath}
\usepackage{amssymb}
\usepackage{mathtools}
\usepackage[hidelinks]{hyperref}
\usepackage{multicol}
\usepackage{makecell}
\usepackage{pgfplots}
\usepackage[]{ulem}
\usepackage{comment}
\usepackage{todonotes}
\usepackage{graphicx}
\usepackage{authblk}
\usepackage{booktabs}
\usepackage[caption = false]{subfig}
\usepackage{stackengine}
\usepackage{algpseudocode, algorithm, algorithmicx}
\usepackage{xr}
\DeclarePairedDelimiter\abs{\lvert}{\rvert}%
\DeclarePairedDelimiter\norm{\lVert}{\rVert}%
\makeatletter
\let\oldabs\abs
\def\abs{\@ifstar{\oldabs}{\oldabs*}}
\let\oldnorm\norm
\def\norm{\@ifstar{\oldnorm}{\oldnorm*}}
\makeatother

\newcommand*\Let[2]{\State #1 $\gets$ #2}
\newcommand*\Compute[2]{\State #1 for #2}

\begin{document}
\normalem 
\title{Regularized Bidimensional Estimation of the Hazard Rate}

\author[1]{Vivien Goepp}
\author[1]{Jean-Christophe Thalabard}
\author[2]{Gr\'egory Nuel}
\author[1]{Olivier Bouaziz}
\affil[1]{MAP5 {\it(CNRS UMR 8145, 45, rue des Saints-P\`eres, 75006 Paris)}}
\affil[2]{LPSM {\it{(CNRS UMR 8001, 4, Place Jussieu, 75005 Paris)}}}
\date{January 2019}
\setcounter{Maxaffil}{0}
\renewcommand\Affilfont{\itshape\small}
\maketitle

\begin{abstract}
In epidemiological or demographic studies, with variable age at onset, a typical quantity of interest is the incidence of a disease (for example the cancer incidence).
In these studies, the individuals are usually highly heterogeneous in terms of dates of birth (the cohort) and with respect to the calendar time (the period) and appropriate estimation methods are needed.
In this article a new estimation method is presented which extends classical age-period-cohort analysis by allowing interactions between age, period and cohort effects.
\textcolor{black}{This paper introduces a bidimensional regularized estimate of the hazard rate where a penalty is introduced on the likelihood of the model.
This penalty can be designed either to smooth the hazard rate or to enforce consecutive values of the hazard to be equal, leading to a parsimonious representation of the hazard rate.
In the latter case, we make use of an iterative penalized likelihood scheme to approximate the L$_0$ norm, which makes the computation tractable.}
The method is evaluated on simulated data and applied on breast cancer survival data from the SEER program.

\end{abstract}

\section*{Introduction}

In epidemiological or demographic studies, with variable age at onset, a typical quantity of interest is the incidence or the hazard rate of a disease (for example the cancer incidence).
In these studies, individuals are recruited and followed-up during a long period of time, usually from birth.
The data are then reported either in the form of registers, which contain the number of observed cases and the number of individuals at risk to contract the disease, or in the form of the observed time for each individual.
These types of studies are of great interest for the epidemiologist, especially when the event of interest will tend to occur at late ages, such as in cancer studies.
However, these data are usually highly heterogeneous in terms of dates of birth and with respect to the calendar time.
In such cases, it is therefore very important to take into account the variability of the age, the cohort (date of birth) and the period (the calendar time) in the hazard rate estimation.
This is usually done using age-period-cohort estimation methods \citep[see][ and citations therein]{Yang2013AgePeriodCohortAnalysisNew}.

In age-period-cohort analysis, the effects of age, period and cohort are fit as factor variables in a regression model where the output is the logarithm of the hazard rate.
However, this induces an identifiability problem due to the relationship: period $=$ age $+$ cohort.
There have been several solutions proposed to this problem.
\cite{Osmond1982AgePeriodCohort} proposed to compute each submodel (age-cohort, age-period, and period-cohort) and use a weighting procedure to combine the three models.
Different constraints have also been proposed to make the age-period-cohort model identifiable.
However, as noticed by \citet[p 162]{Heuer1997ModelingTimeTrends}, the obtained estimates highly depend on the choice of the constraints.
\cite{Holford1983EstimationAgePeriod} proposed to directly estimate the linear trends of each effect.
This procedure leads to results that are difficult to interpret.
See \cite{Carstensen2007AgePeriodCohort} for a detailed discussion of the identifiability problem of the age-period-cohort model.
More recently,  \cite{Kuang2008IdentificationAgePeriodCohortModel} proposed to estimate the second order derivatives of the three effects.
This model is implemented in the package \texttt{apc} \cite{Nielsen2015apcPackageAgePeriodCohort}.
Finally, \cite{Carstensen2007AgePeriodCohort} proposed to first fit one submodel (say age-cohort) and then to fit the period effect over the residuals of the  first model. 
This model is implemented in the \texttt{R} package \texttt{Epi} \citep{CarstensenEpiPackageStatistical2017}, \cite{Plummer2011LexisClassEpidemiological}.

All these approaches can be viewed as parametric models, where the parameters are the age, period, and cohort vector parameters.
As such they are also restrictive because they do not allow for interactions between the three effects, that is they assume that one effect does not depend on the other effect's value. 
A different approach consists in considering the hazard rate as a function of age and either period or cohort and to estimate this bi-dimensional function in a non-parametric setting.
No specific structure of the hazard rate is assumed.
\textcolor{black}{However, for moderate sample sizes, non-parametric approaches such as the maximum likelihood estimator (MLE) are prone to overparametrization. As a matter of fact the MLE can only be used if a bi-dimensional grid (e.g. of cohort and age intervals) is provided. Without an appropriate method, this grid needs to be arbitrarily chosen. If the number of intervals is too large, the MLE will display a large variance. On the other hand, a too small number of intervals will result in a bias if those intervals are not optimally chosen.}
Consequently, regularized methods have been proposed in order to avoid overfitting in this non-parametric context.
A kernel-type estimator was proposed by \cite{Beran1981NonparametricRegressionRandomlya} and \cite{McKeague1990IdentifyingNonlinearCovariate} where the cumulative hazard is smoothed using a kernel function.
See \cite{Keiding1990StatisticalinferenceLexis} for a thorough discussion of methods for hazard inference in age-period-cohort analysis.
More recently, \cite{Currie2009SmoothingAgePeriodCohortModels} proposed a spline estimation procedure to infer the hazard rate as a function of two variables.
The authors use a generalized linear model using B-splines and overfitting is dealt with using a penalization over the differences of adjacent splines' coefficients.

In this article, we propose a new non-parametric method for bi-dimensional hazard rate estimation.
As the previous non-parametric approaches, this model considers the estimation of the hazard rate with respect to two variables, i.e. either age-cohort, age-period, or period-cohort, without assuming any specific structure on the hazard rate. 
Inference is made in two dimensions, but through the linear relationship period $=$ age $+$ cohort, the hazard rate can be represented as a function of any two of the three variables.
\textcolor{black}{Finally, in order to take into account the issue of overfitting, we use a sparsity-inducing penalized likelihood method called adaptive ridge. This iterative method is an approximation of the L$_{0}$ norm penalty which makes the computation tractable. We note that the L$_0$ ``norm'', defined by $\|\bm{u}\|_{0} = \#\{j|u_j \neq 0\}$ is not a proper norm but we nevertheless use the term ``L$_0$ norm'' hereafter following the notations of e.g. \cite{Candes2008EnhancingSparsityReweighted}. The method was first introduced by \cite{Chartrand2008Iterativelyreweightedalgorithms} in the context of sparse sensing and applied by \cite{Rippe2012VisualizationGenomicChanges} and \cite{Frommlet2016AdaptiveRidgeProcedure} in the context of linear regression. It has been used in the context of piecewise constant hazard rate estimation by  \cite{Bouaziz2017L0RegularizationEstimation}. The present work makes use of this method to perform a segmentation of the hazard rate into constant areas.}
The novelty of this method lies in the parsimonious representation of the bi-dimensional hazard rate into segmented areas.
In particular, the method can efficiently exhibit cohort, age or period effects, that is, specific changes of the hazard rate due to the date of birth, the age or the calendar time.
The penalized likelihood framework used here can also be used to estimate the L$_2$ norm penalization, which will induce a smoothed estimate of the hazard in a similar way as the aforementioned non-parametric methods.

Our model is introduced in Section~\ref{section:model}.
The regularization method is then presented in Section~\ref{section:regul}.
In Section~\ref{section:penalty}, the selection of the penalty parameter is discussed.
Finally, the performance of our model is assessed through a simulation study in Section~\ref{section:simu} and illustrations on the SEER cancer dataset is provided in Section~\ref{section:realdata}.

\begin{figure}[H]
\centering
\begin{tabular}{cc}
\begin{tikzpicture}[scale = 0.8, xscale = 0.07, yscale = 0.07, domain = 0.140:60, samples = 800, every node/.style = {scale = 0.5}]
\tikzstyle{instruct} = [rectangle, draw, fill = yellow!50]
\draw[->] (0, 0) -- (110, 0) node[right] {period};
\draw[->] (0, 0) -- (0, 110) node[left] {age};
\foreach \i in {1900, 1910, 1920, ..., 2000} {
  \draw (\i - 1900, 0) -- (\i - 1900, -2) node[below] {\small $\i$};
  \draw [lightgray, dashed](\i - 1900, 0) -- (\i - 1900, 110);
  \draw [lightgray, dashed](\i - 1900, 0)-- (100, 100 - \i + 1900);
}
\foreach \i in {0, 10, ..., 100} {
  \draw (0, \i) -- (-2, \i) node[left] {$\i$};
  \draw [lightgray, dashed](0, \i) -- (110, \i);
  \draw [lightgray, dashed](0, \i) -- (100 - \i, 100);
}
\draw (0, 0) -- (110, 0);
\draw (0, 0) -- (0, 110);
\draw (23, 0) node {$\circ$} -- (74, 51) node {$\times$};
\draw (55, 0) node {$\circ$} -- (90, 37) node {$\times$};
\fill [opacity = 0.3] (30, 0) -- (100, 70) -- (100, 80) -- (20, 0) -- cycle;    
\fill [opacity = 0.2] (50, 0) -- (50, 105) -- (60, 105) -- (60, 0) -- cycle;
\fill [opacity = 0.2] (0, 30) -- (100, 30) -- (100, 40) -- (0, 40) --cycle;
\node[rectangle, draw, fill = black!20] (plus) at (95, 69) {cohort};
\node[rectangle, draw, fill = black!20] (plus) at (55, 100) {period};
\node[rectangle, draw, fill = black!20] (plus) at (100, 35) {\phantom{ii}age\phantom{ii}};
\end{tikzpicture}
&
\begin{tikzpicture}[scale = 0.8, xscale = 0.07, yscale = 0.07, domain = 0.140:60, samples = 800, every node/.style = {scale = 0.5}]
\draw[->] (0,0) -- (110,0) node[right] {cohort};
\draw[->] (0,0) -- (0,110) node[left] {age};
\foreach \i in {1900, 1910, 1920, ..., 2000} {
  \draw (\i - 1900, 0) -- (\i - 1900, -2) node[below] {\small $\i$};
  \draw [lightgray, dashed](\i - 1900, 0) -- (\i - 1900, 110);
  \draw [lightgray, dashed](\i - 1900, 100)-- (100, \i - 1900);
}
\foreach \i in {0, 10, ..., 100} {
  \draw (0, \i) -- (-2, \i) node[left] {$\i$};
  \draw [lightgray, dashed](0, \i) -- (110, \i);
  \draw [lightgray, dashed](0, \i) -- (\i, 0);
}
\draw (0, 0) -- (110, 0);
\draw (0, 0) -- (0, 110);
\draw (23, 0) node {$\circ$} -- (23, 51) node {$\times$};
\draw (55, 0) node {$\circ$} -- (55, 37) node {$\times$};
\fill [opacity = 0.2] (50, 0) -- (60, 0) -- (0, 60) -- (0, 50) -- cycle;   
\fill [opacity = 0.2] (20, 0) -- (30, 0) -- (30, 100) -- (20, 100) -- cycle;
\fill [opacity = 0.2] (0,30) -- (100, 30) -- (100, 40) -- (0, 40) -- cycle;
\node[rectangle, draw, fill = black!20] (plus) at (25, 100) {cohort};
\node[rectangle, draw, fill = black!20] (plus) at (7, 49) {period};
\node[rectangle, draw, fill = black!20] (plus) at (100, 35) {\phantom{ii}age\phantom{ii}};
\end{tikzpicture} \\
(a) Lexis diagram: Age-Period diagram & (b) Age-Cohort diagram \\
\end{tabular}
\caption{Diagrams representing the lives of individuals: in the age-period plane (a) -- called Lexis diagram -- and in the age-cohort plane (b). Solid lines represent lives of individuals until occurrence of the event of interest. The same age, cohort, and period intervals are displayed in gray.}
\label{fig:diagram}
\end{figure}
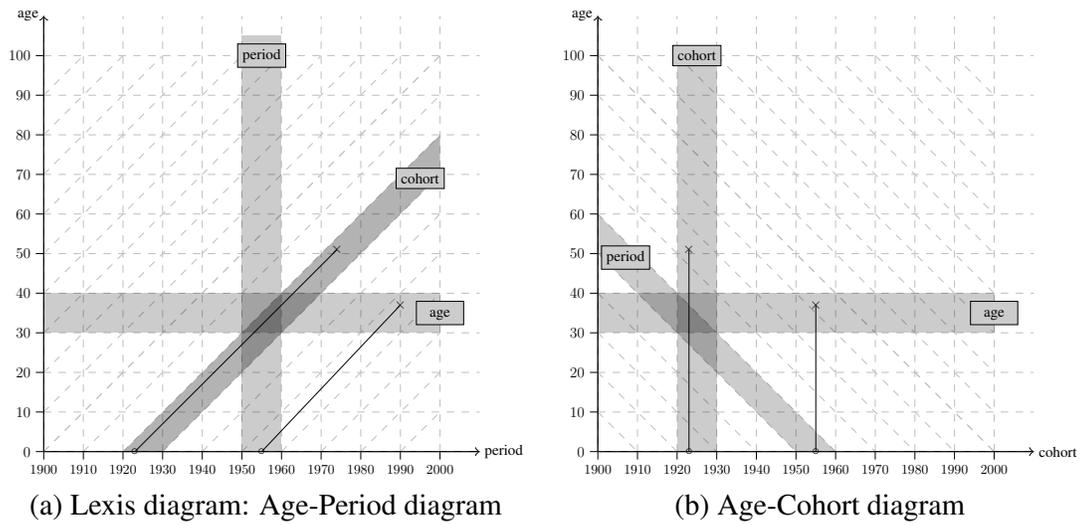

\textcolor{black}{\section{Fused Regularized Estimation}%
\label{section:model}
\subsection{Modelization}}
In the age-period-cohort setting, the date of birth (the cohort) $U$ of each individual is available and the variable of interest is a time-to-event variable of this individual denoted $T$.
The data are subject to right-censoring and are represented as tabulated data over the $J$ cohort intervals and the $K$ age intervals $[c_0, c_1), [c_1, c_2), \dots, [c_{J-1}, c_J)$ and $[d_0, d_1), [d_1, d_2), \dots, [d_{K-1}, d_K)$ respectively, with the convention $c_0=d_0=0$ and \textcolor{black}{$c_J=d_K = \infty$}.
On a sample of $n$ individuals, the available data can then be rewritten in terms of the exhaustive statistics $\bm{O}=(O_{1,1},\ldots,O_{J,K})$, $\bm{R}=(R_{1,1},\ldots,R_{J,K})$, where for $j=1,\ldots,J$, $k=1,\ldots,K$, $O_{j,k}$ represents the number of observed events that occurred in the $j$-th cohort interval $[c_{j-1}, c_j)$ and $k$-th age interval $[d_{k-1}, d_k)$ and $R_{j,k}$ represents the total time individuals were at risk in this $j$-th cohort and $k$-th age interval.
In the case of register data, the discretization $\left( c_{j} \right), \left( d_{k} \right)$ is imposed by the data and the available data is directly $\bm{R}$ and $\bm{O}$, which are often called the \emph{cases} and \emph{person-years}, respectively.
See for instance \cite{Carstensen2007AgePeriodCohort} for an example of such data.
The aim is to estimate the hazard rate, defined as:  
\begin{equation*}
\lambda(t| u) = \lim_{dt \to 0} \frac{1}{dt} \text{P}(t < T < t + dt | T > t, U = u).
\end{equation*}
In the age-cohort setting \textcolor{black}{$\lambda(t| u)$} is assumed to be piecewise constant:
\begin{equation*}
\lambda(t| u) = \sum_{j=1}^J \sum_{k = 1}^K \lambda_{j, k} \mathds{1}_{[c_{j-1}, c_j) \times [d_{k-1}, d_k)}(t,u),
\end{equation*}
and inference is made over the $J\times K$ dimension parameter $\bm{\lambda}=(\lambda_{1, 1},\ldots,\lambda_{J, K})$.
Note that the hazard can be equivalently defined as a function of age and period or as a function of period and cohort where the period is defined as the calendar time, that is: period $=$ cohort $+$ age.
For illustration, the change of coordinates between the age-period and age-cohort diagrams is represented in Figure \ref{fig:diagram}.
In our models, the hazard will be considered as a function of solely age and cohort since the influence of any of the two elements of age, period or cohort can be retrieved using this reparametrization.

\textcolor{black}{\subsection{Penalized Likelihood}}

Following \citet[p. 224]{Aalen2008SurvivalEventHistory} the negative log-likelihood takes the form
\begin{equation}
		\ell_n(\bm{\lambda}) = \sum_{j = 1} ^ J \sum_{k = 1} ^ K \{ \lambda_{j, k} R_{j,k} - O_{j,k} \log \left(\lambda_{j, k}\right)\}.
\label{nllh}
\end{equation}
The authors also noticed that this log-likelihood is equivalent to a log-likelihood arising from a Poisson model. However, note that no distribution assumptions are made on the data and in particular the $O_{j,k}$ are not assumed to be Poisson distributed \citep[see][for a discussion on the ``Poisson'' model]{Carstensen2007AgePeriodCohort}.
Minimizing $\ell_n$ yields an explicit maximum likelihood estimate $\widehat{\lambda}^{\text {mle}}_{j, k} = O_{j,k} / R_{j,k}$. However, for moderate sample sizes this estimator is overfitted, especially in places of the age-cohort plane where few events are recorded. 
To remedy this problem we propose in the following to penalize the differences between adjacent values of the hazard in the log-likelihood.

For computation convenience, we first reparametrize the model: $\eta_{j,k} = \log \lambda_{j,k}$, for $1 \leq j \leq J$ and $1 \leq k \leq K$.
The goal of this work is to estimate the minimizer of the function
\begin{align}
\ell_n(\bm{\eta}) + \frac{\kappa}{2} \sum_{j=1}^{J - 1} \sum_{k = 1} ^ {K}  \norm{\eta_{j + 1, k} - \eta_{j, k}}_0 + \frac{\kappa}{2} \sum_{j=1}^{J} \sum_{k = 1} ^ {K - 1} \norm{\eta_{j, k + 1} - \eta_{j, k}}_0,
\label{eq:rem_nll_l0}
\end{align}
where $\ell_n(\bm{\eta})$ was defined in Equation~\eqref{nllh} and $\kappa$ is a penalty constant used as a tuning parameter.
In the previous equation, the L$_0$ norm penalty  over the differences of adjacent parameter values yields a piecewise constant estimate of the hazard rate.
However this function is not tractable to minimize due to the L$_0$ norm.
Hence we use the adaptive ridge, which, as is explained in Section~\ref{section:regul}, can be seen as an approximate solution to this problem.

Let us define the weighted L$_2$ penalized model:
\begin{align}
\ell_n^\kappa(\bm{\eta}, \bm{v}, \bm{w}) = \ell_n(\bm{\eta}) + \frac{\kappa}{2} \sum_{j=1}^{J - 1} \sum_{k = 1} ^ {K}  v_{j, k} \left( \eta_{j + 1, k} - \eta_{j, k} \right) ^ 2 + \frac{\kappa}{2} \sum_{j=1}^{J} \sum_{k = 1} ^ {K - 1} w_{j, k} \left( \eta_{j, k + 1} - \eta_{j, k} \right) ^ 2,
\label{eq:rem_nll}
\end{align}
where $\bm{v} = (v_{1,1},\ldots,v_{J-1,K})$ and $\bm{w} = (w_{1,1},\ldots,w_{J,K-1})$ are constant positive weights of respective dimensions $(J-1)K$ and $J(K - 1)$.
Note that the case $\kappa = 0$ corresponds to the maximum likelihood estimation and the case $\kappa = \infty$ corresponds to a hazard uniformly constant over the age and cohort intervals.
The parameter $\kappa$ needs to be chosen in an appropriate way in order to obtain a compromise between these two extreme situations. This is addressed in Section~\ref{section:penalty}.

This model does not attempt to estimate the age, period and cohort effect as parameter vectors.
Instead, it performs a regularized estimation of $\bm \lambda$ that has no age-period-cohort-type structure.

In the next section, we introduce an algorithm to minimizing Equation~\eqref{eq:rem_nll}, which will be used for estimating both L$_2$ and L$_0$ penalties. We then introduce the estimation procedures for both fused L$_2$ and L$_0$ penalties.

\textcolor{black}{\section{Numerical Optimization}%
\label{section:regul}}
In this section, we first introduce the weighted L$_2$ penalized negative log-likelihood and derive how to minimize it.
Then, two different expressions of the weights $\bm{v}$ and $\bm{w}$ are proposed which correspond to two different types of regularization of the hazard rate. The first one implements the adaptive ridge and yields a piecewise constant estimate. The second one uses constant weights and yields a smooth estimate.

\textcolor{black}{\subsection{Fused L$_{2}$ Penalty Estimate}}

Minimization of $\ell_n^\kappa$ is performed using the Newton-Raphson method (see Algorithm~\ref{alg:nr}). Let $\bm{U}_n^\kappa(\bm{\eta}, \bm{v}, \bm{w}) = \partial \ell_n^\kappa/\partial\bm{\eta}$ be the gradient of the penalized negative log-likelihood and $\bm{I}_n^\kappa(\bm{\eta}, \bm{v}, \bm{w}) = \partial \bm{U}_n^\kappa(\bm{\eta}, \bm{v}, \bm{w}) / \partial \bm{\eta}^T$ be its Hessian matrix. 

For $1\leq j,j' \leq J$ and $1\leq k,k' \leq K$, simple algebra yields
\begin{equation*}
\frac{\partial \ell_n(\bm{\eta})}{\partial \eta_{j, k}} = \exp \left(\eta_{j, k}\right) R_{j,k} - O_{j, k}, \quad \quad \frac{\partial ^2 \ell_n(\bm{\eta})}{\partial \eta_{j', k'} \partial \eta_{j, k}} = \mathds{1}_{j=j', k=k'} \exp \left(\eta_{j, k}\right)R_{j,k}, \quad \text{and}
\end{equation*}
\begin{align*}
\frac{\partial \ell_n ^\kappa}{\partial \eta_{j,k}}(\bm{\eta}) &= \frac{\partial \ell_n(\bm{\eta})}{\partial \eta_{j,k}} + \kappa \left[ - v_{j, k} \left( \eta_{j+1, k} - \eta_{j,k}\right) + v_{j-1, k} \left(\eta_{j, k} - \eta_{j - 1,k} \right)\right] \\
& \quad \quad \quad \quad \phantom{i.} + \kappa \left[ - w_{j, k} \left( \eta_{j, k+1} - \eta_{j,k}\right) + w_{j, k - 1} \left(\eta_{j, k} - \eta_{j, k - 1} \right)\right],
\end{align*}
\begin{align*}
\frac{\partial ^2\ell_n ^ \kappa(\bm{\eta})}{\partial \eta_{j', k'} \partial \eta_{j, k}} & = \frac{\partial ^2\ell_n(\bm{\eta})}{\partial \eta_{j', k'} \partial \eta_{j, k}} + \kappa \left[ \mathds{1}_{j = j', k = k'} \left( v_{j', k'} + v_{j' - 1, k'} + w_{j', k'} + w_{j', k' - 1} \right)\right. \\
& \phantom{MMMMMMMMMM} - v_{j', k'} \mathds{1}_{j = j' + 1, k = k'} -  v_{j' - 1, k'} \mathds{1}_{j = j' - 1, k = k'} \\
& \phantom{MMMMMMMMMM}\left. - w_{j', k'} \mathds{1}_{j = j', k = k' + 1} -  w_{j', k' - 1} \mathds{1}_{j = j', k = k' - 1} \right].
\end{align*}
From the last equation, the Hessian matrix can be written
\begin{equation*}
\bm{I}_n^\kappa(\bm{\eta}, \bm{v}, \bm{w}) = \frac{\partial^2 \ell_n(\bm{\eta})}{\partial \bm{\eta} \partial \bm{\eta}^T} + \kappa  \bm{B}(\bm{\eta}),
\label{hessian}
\end{equation*}
where $\bm{B}(\bm{\eta})$ is a band matrix of bandwidth equal to $\min(J, K) - 1$. \textcolor{black}{Thus the Hessian matrix is also a band matrix of bandwidth $\min(J, K) - 1$}. Using Cholesky decomposition, the computation of $\bm{I}_n^\kappa(\bm{\eta}, \bm{v}, \bm{w}) ^{-1} \bm{U}_n^\kappa(\bm{\eta}, \bm{v}, \bm{w})$ has a $\mathcal{O}(\min(J, K) J K)$ complexity instead of $\mathcal{O}(J^3K^3)$.

\begin{algorithm}
  \caption{Newton-Raphson Procedure with Constant Weights
    \label{alg:nr}}
  \begin{algorithmic}[1]
    \Statex
    \Function{Newton-Raphson}{$\bm{O}, \bm{R},\kappa, \bm{v}, \bm{w}$}
      \Let{$\bm{\eta}$}{$\bm{0}$}
      \While{not converge}
        \Let{$\bm{\eta}^{\text{new}}$}{$\bm{\eta} - \bm{I}_n^\kappa (\bm{\eta}, \bm{v}, \bm{w})^ {-1} \bm{U}_n^\kappa(\bm{\eta}, \bm{v}, \bm{w})$}
        \Let{$\bm{\eta}$}{$\bm{\eta}^{\text{new}}$}
      \EndWhile
      \State \Return{$\bm{\eta}$}
    \EndFunction
  \end{algorithmic}
\end{algorithm}

A ridge-type penalization is performed when setting $\bm{v} = \bm{w} = \bm{1}$ in Equation \eqref{eq:rem_nll_l0}.
In this case the penalization corresponds to the square of the first-order differences of  \textcolor{black}{$\bm{\eta}$}, which yields a smooth estimator of the hazard rate.
{\color{black}This estimate is obtained directly from Algorithm \ref{alg:nr}.}

We make a note that Equation \ref{eq:rem_nll} allows for some flexibility in the regularization.
Indeed, one could set different values to $\bm{v}$ and $\bm{w}$ to manually tune the importance of the regularization between different regions of the plane and between the two variables.

Finally, note that this method yields an estimate similar to the spline method of \cite{Ogata1988LikelihoodAnalysisSpatial}, who penalizes over the second-order differences instead of the first-order differences.
This means that for arbitrarily large values of the penalty constant, the regularized hazard will be a constant function instead of a linear function.


{\color{black}
\subsection{Fused Adaptive Ridge Estimate}
{\color{black}
In this section, we derive a computationally tractable procedure to minimize Equation \eqref{eq:rem_nll_l0}.
We make use of the adaptive ridge, which minimizes a non-convex penalty by iteratively minimizing approximations of the penalty.
The adaptive ridge can be used to approximate any L$_q$ penalty ($0<q<1$) and it extends to the case $q = 0$, the latter case corresponding to the logarithmic penalty in lieu of the L$_0$ penalty.
This procedure is still called ``L$_0$ adaptive ridge'' since, as explained by \cite{Candes2008EnhancingSparsityReweighted}, the logarithmic penalty is a good approximation of the L$_0$ penalty.
The adaptive ridge iteratively solves L$_2$ penalty problems (hence its name), and is thus simple to implement.

As pointed out by a reviewer, another iterative penalized method \citep{Foucart2009Sparsestsolutionsunderdetermined} could have been used which iterately solves L$_1$ penalty problems to approximate the L$_0$ penalty.
As explained by \cite{Wipf2010Iterative}, these two methods are very similar in that they both minimize a logarithm penalized problem using two different approximations.
In Section~1 of Supplementary Material, we make the link between the two methods explicit and show that they belong to the same class of optimization schemes. We also refer there to related works using either of the two methods.
}

We implement the adaptive ridge procedure by minimizing Equation \eqref{eq:rem_nll} with the weights adapted iteratively.
We iterate between updating 
\begin{equation*}
    \bm{\eta}^\text{new} \gets \arg\min_{\bm{\eta}} \ell_n^\kappa(\bm{\eta}, \bm{v}, \bm{w})
\end{equation*}
using Algorithm~\ref{alg:nr} and updating the values of the weights:
\begin{equation*}
\begin{cases}
v_{j, k}^{\text{new}} = \left(\left(\eta_{j + 1, k} ^ {\text{new}} - \eta_{j, k} ^ {\text{new}}\right) ^ 2 + \varepsilon_v ^ 2\right) ^{-1}, \\ 
w_{j, k}^{\text{new}} = \left(\left(\eta_{j, k} ^ {\text{new}} - \eta_{j, k - 1} ^ {\text{new}}\right) ^ 2 + \varepsilon_w ^ 2\right) ^ {-1}, \\
\end{cases}
\label{eq:weights}
\end{equation*}
where $\varepsilon_v$ and $\varepsilon_w$ are constants negligible compared to 1. 
}

We now elaborate on the estimation procedure.
{\color{black} The algorithm is said to converge if the absolute difference of all weighted differences in $\eta_{j,k}$ are below a given threshold (we use $10^{-8}$ in our implementation).
At convergence, $v_{j, k} \left(\eta_{j+1, k} - \eta_{j, k}\right)^2$ will be either very close to $0$ if $|\eta_{j+1, k} - \eta_{j, k}|$ is smaller than $\varepsilon_{v}$ or very close to $1$ if $|\eta_{j+1, k} - \eta_{j, k}|$ is greater than $\varepsilon_{v}$ -- and similarly for $w_{j, k}\left(\eta_{j, k+1} - \eta_{j, k}\right)^2$.
We then set them to $0$ or $1$ using a thresholding, so that values smaller than $0.99$ are set to $0$ and values larger than $0.99$ are set to $1$ \textcolor{black}{(in practice, the value of this threshold has little effect, since at convergence the weighted differences are distant to $0$ or $1$ by $\sim 10^{-7}$)}.

As with other penalized methods and as pointed out in \cite{Frommlet2016AdaptiveRidgeProcedure}, the adaptive ridge penalization scheme induces a shrinkage bias.
Therefore, after segmentation of the $\eta_{i,j}$s, the hazard rate is estimated on each constant area using the unpenalized maximum likelihood estimator.
These constant areas are defined as connected components of a graph.
We first create the graph whose vertices are the $JK$ age-cohort rectangles and whose edges are the connections between adjacent cells that have differences equal to $0$. Then, each connected component of this graph is a different area over which the hazard has been estimated to be constant.}
The extraction of connected components from the graph is done using the package \texttt{igraph} \citep{Csardi2006igraphSoftwarePackage}. 
The log-hazard $\eta^{(r)}$ of the $r$-th constant area is such that $\forall [c_{j-1}, c_j) \times [d_{k - 1}, d_k) \in r,  \eta_{j,k} = \eta^{(r)}$.
\textcolor{black}{The values of $\eta^{(r)}$ are then estimated in a second step, using unpenalized maximum likelihood estimation: $\widehat{\eta}^{(r)} = \log \left(O^{(r)} / R^{(r)}\right)$ where $O^{(r)}$ is the number of events in the $r$-th constant area and $R^{(r)}$ is the time at risk in the $r$-th constant area.}

This algorithmic procedure is summarized in Algorithm \ref{alg:aridge}.
In practice, the stopping criterion for the adaptive ridge algorithm is when the absolute difference between successive values of the weighted differences is smaller than a predefined value -- we use $10^{-8}$ in our implementation.
Moreover, following \cite{Frommlet2016AdaptiveRidgeProcedure}, we have set $\varepsilon_v = \varepsilon_w = 10 ^ {-5}$.

\begin{algorithm}[h]
  \caption{Adaptive Ridge Procedure
    \label{alg:aridge}}
  \begin{algorithmic}[1]
    \Statex
    \Function{Adaptive-Ridge}{$\bm{O}, \bm{R}, \kappa$}
      \Let{$\bm{\eta}$}{$\bm{0}$}
      \Let{$\bm{v}$}{$\bm{1}$}
      \Let{$\bm{w}$}{$\bm{1}$}
      \While{not converge}
        \Let{$\bm{\eta}^{\text{new}}$}{\textsc{Newton-Raphson}$(\bm{O}, \bm{R}, \kappa, \bm{v}, \bm{w})$}
        
\Let{$v_{j, k}^{\text{new}}$}{$\left(\left(\eta_{j + 1, k} ^ {\text{new}} - \eta_{j, k} ^ {\text{new}}\right) ^ 2 + \varepsilon_v ^ 2\right) ^ {-1}$}
\Let{$w_{j, k}^{\text{new}}$}{$\left(\left(\eta_{j, k} ^ {\text{new}} - \eta_{j, k - 1} ^ {\text{new}}\right) ^ 2 + \varepsilon_w ^ 2 \right) ^ {-1}$}
        \Let{$\bm{\eta}$}{$\bm{\eta}^{\text{new}}$}
      \EndWhile
      
      \Compute{{\bf Compute} $(O^{\text{new}}, R^{\text{new}})$}{selected $(\bm{\eta}, \bm{v}^{\text{new}}, \bm{w}^{\text{new}})$}
      \Let{$\bm{\eta}^{\text{new}}$}{$\log \left( \bm{O}^{\text{new}}/\bm{R}^{\text{new}}\right)$}
      \State \Return{$\bm{\eta}^{\text{new}}$}
    \EndFunction
  \end{algorithmic}
\end{algorithm}

\section{Choice of the Penalty Parameter $\kappa$}%
\label{section:penalty}
In practice, the hazard rate needs to be estimated for a set of penalty constants and the choice of $\kappa$ is determined as the penalty that provides the best compromise between model fit and reduced variability of the hazard rate estimate. For the L$_0$ regularization model, different values of the penalty constant lead to different segmentations of the $\eta_{j,k}$. As a consequence, the problem of choosing the optimal penalty constant can be rephrased as the problem of choosing the optimal model among a set of models $\mathcal{M}_1, \dots, \mathcal{M}_M$, where each of these models corresponds to a different segmentation of the $\eta_{j,k}$ and $M$ is the maximum number of different models. In this section we propose different methods to select the optimal model. Comparison of the efficiency of the different methods will be analyzed in Section \ref{section:simu} on simulated data.

We recall that $\bm{R}$ and $\bm{O}$ are the exhaustive statistics and $\bm{\eta}$ is the parameter to be estimated in our two models. Bayesian criteria attempt to maximize the posterior probability $\text{P}(\mathcal M_m | \bm{R}, \bm{O}) \propto \text{P}(\bm{R}, \bm{O}|\mathcal M_m) \pi(\mathcal M_m)$, where $\text{P}(\bm R, \bm O | \mathcal{M}_m)$ is the integrated likelihood and $\pi\left( \mathcal{M}_m \right)$ is the prior distribution on the model. This problem is equivalent to minimizing $-2 \log \text{P}(\mathcal{M}_m | \bm{R}, \bm{O})$. By integration 
\begin{equation*}
\text{P}(\bm{R}, \bm{O}|\mathcal M_m) = \int_{\bm{\eta}} \text{P}(\bm{R}, \bm{O}|\mathcal M_m, \bm \eta) \pi(\bm \eta) d\bm \eta, 
\end{equation*}
where $\text{P}\left(\bm R, \bm O | \mathcal{M}_m, \bm \eta\right)$ is the likelihood and $\pi(\bm \eta)$ is the prior distribution of the parameter, which is taken constant in the following.
Thus Bayesian criteria are defined as
\begin{equation*}
- 2 \log \left( \text{P} (\mathcal{M}_m | \bm R, \bm O\right) = 2 \ell_n(\widehat{\bm{\eta}}_m) + q_m \log n - 2 \log \pi(\mathcal M_m) + \mathcal{O}_{\text{P}}(1),
\end{equation*}
where $q_m$ is the dimension of the model $\mathcal{M}_m$ i.e., the number of constant areas selected by the adaptive ridge algorithm.

The BIC \citep{Schwarz1978EstimatingDimensionModel} corresponds to the Bayesian criterion obtained when one neglects the term $\pi(\mathcal M_m)$, which is equivalent to having a uniform prior on the model:

\begin{equation}
\text{BIC}(m) = 2 \ell_n (\widehat{\bm{\eta}}_m) + q_m \log n.
\label{bic}
\end{equation}

As explained by \cite{Zak-Szatkowska2011ModifiedVersionsBayesian}, a uniform prior on the model is equivalent to a binomial prior on the model dimension $\mathcal{B}(JK, 1/2)$.
When the true model's dimension is much smaller than the maximum possible dimension $JK$, the BIC tends to give too much importance to models of dimensions around $JK / 2$, which will result in underpenalized estimators.
To this effect, \cite{Chen2008ExtendedBayesianInformation} have developed an extended Bayesian information criterion called EBIC$_{0}$ (or EBIC for short).
One can write $\pi(\mathcal{M}_m) = \text{P}(\mathcal{M}_m | \mathcal{M}_m \in \mathcal{M}_{[q_m]}) \text{P}(\mathcal{M}_m \in \mathcal{M}_{[q_m]})$ where $\mathcal{M}_{[q_m]}$ is the set of models of dimension $q_m$.
The EBIC$_{0}$ criterion is defined by setting $\text{P}(\mathcal{M}_m | \mathcal{M}_m \in \mathcal{M}_{[q_m]}) = 1 / {J K \choose q_m}$ and $\text{P}(\mathcal{M}_m \in \mathcal{M}_{[q_m]}) = 1$.
Thus 

\begin{equation*}
\pi(\textcolor{black}{\mathcal{M}}_m) = {JK \choose q_m}
\end{equation*} 
and
\begin{equation}
		\text{EBIC}_{0}\left( m \right) = 2 \ell_n (\widehat{\bm{\eta}}_m) + q_m \log n + 2 \log {JK \choose q_m}.
\label{ebic}
\end{equation}

Note that the EBIC$_{0}$ assigns the same \emph{a priori} probability to all models of same dimension.
Therefore, when the true model's dimension is not close to $JK /2$ the EBIC$_{0}$ will be able to select this model more easily.
Namely, when the true model's dimension is very small the EBIC$_{0}$ will tend to choose very sparse models.

The last criterion that will be used is the Akaike Information Criterion \citep{Akaike1998InformationTheoryExtension}, or AIC, defined as $\text{AIC}(m) = 2 \ell_n (\widehat{\bm{\eta}}_m) + 2 q_m$.
This criterion is known for performing better than the BIC in terms of mean squared error, however the BIC will tend to select sparser models than the AIC.

Note that Bayesian criteria and the AIC can only be used for the L$_0$ regularized estimation only, since the L$_2$ model does not perform a model selection.
An alternative to performing model selection is to use the K-fold cross validation. With this method, the data are split at random into $L$ parts. 
The estimated parameter obtained when the $l$-th part is left out is noted $\widehat{\eta}^{-l}(\kappa)$ and the cross-validated score is defined as 
\begin{equation*}
\text{CV}(\kappa) = \sum_{l = 1}^L \ell_n^{\kappa, l}(\widehat{\bm \eta}^{-l}),
\end{equation*}
where $\ell_n^{\kappa, l}$ is the negative log-likelihood evaluated on the $l$-th part of the data.
The optimal penalty constant is obtained by minimizing $\text{CV}(\kappa)$ with respect to $\kappa$.
The L-fold cross validation method can be used for both the L$_0$ regularized estimation and the L$_2$ regularized estimation.
However, this method is numerically time consuming as the estimator has to be computed $L$ times while Bayesian criteria or the AIC provide direct methods to perform model selection from the original estimator.
In the simulation studies and data analysis, we set $L = 10$.

\section{Simulation Study}
\label{section:simu}

\subsection{Simulation Designs}
\label{section:simu_design}

\begin{figure}
\vspace{-0.4 cm}
\begin{center}
$
\begin{array}{cc}
\includegraphics[width = 65mm, height = 65mm]{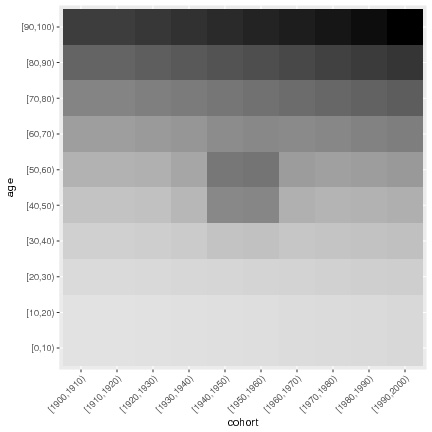} & 
\includegraphics[width = 65mm, height = 65mm]{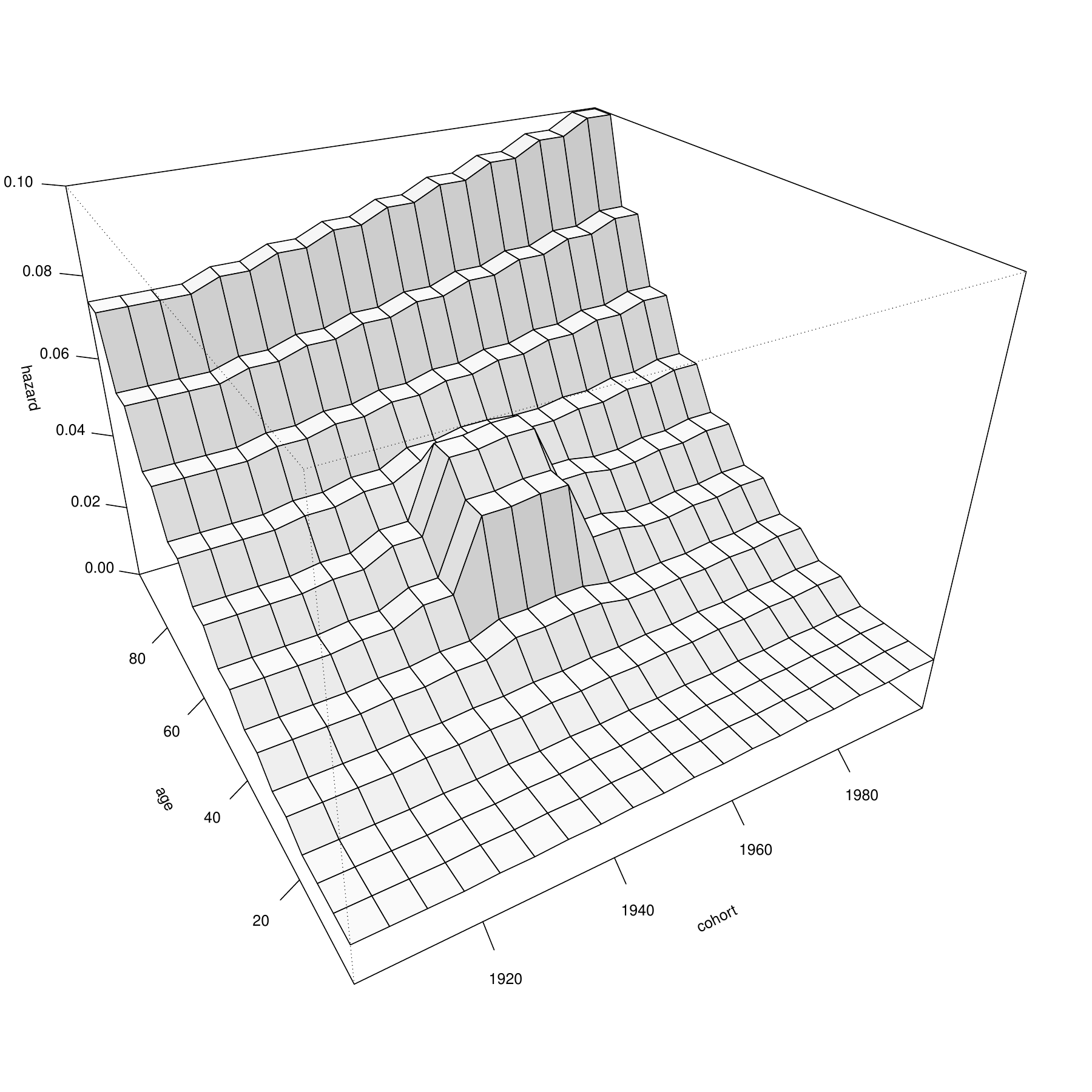} \\
\text{(a) Smooth true hazard -- heatmap} & \text{(b) Smooth true hazard -- perspective} \\
\includegraphics[width = 65mm, height = 65mm]{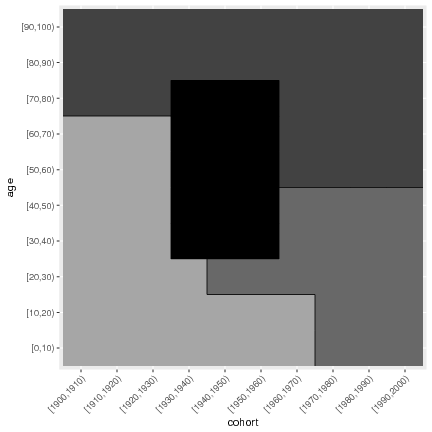} & 
\includegraphics[width = 65mm, height = 65mm]{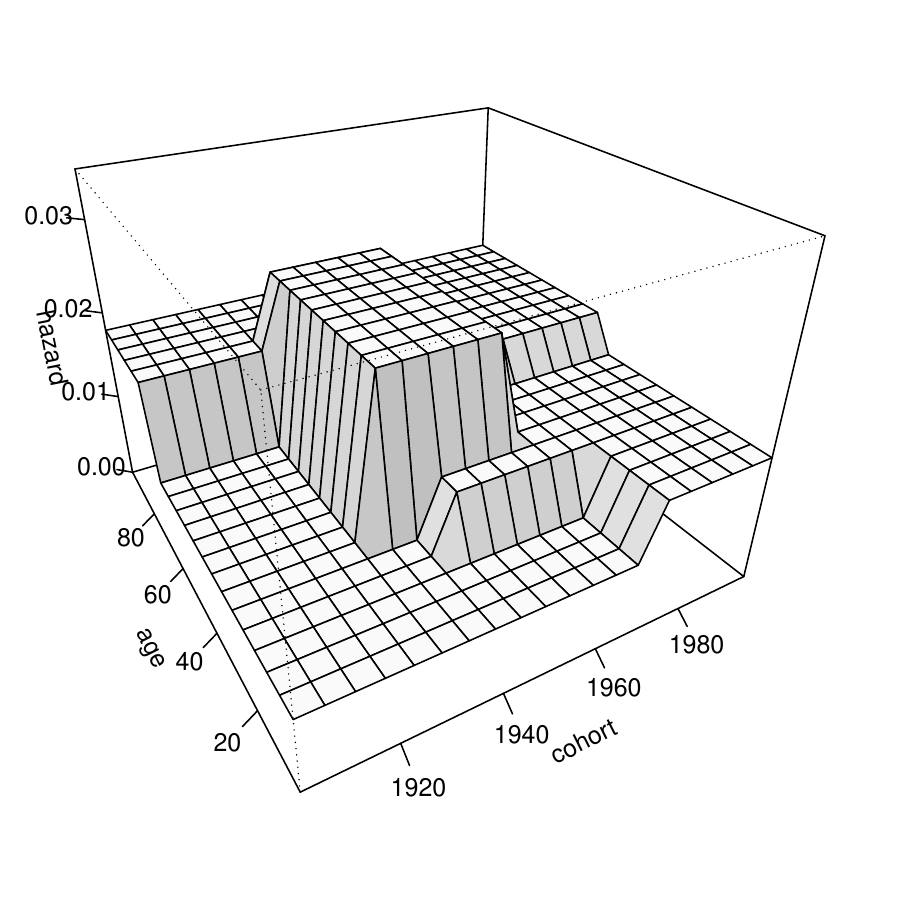}\\
\text{(c) Piecewise constant true hazard -- heatmap} & \text{(d) Piecewise constant true hazard -- perspective}\\
\end{array}$
\end{center}
\caption{True hazard of the two simulation designs: smooth hazard in heatmap (a) and perspective plot (b) and piecewise constant hazard in heatmap (c) and perspective plot (d).}
\label{fig:simu_design}
\end{figure}

In this section, our piecewise estimation method is compared with the \textsc{age-cohort} model and with the L$_2$ penalty estimate.
The different criteria for model selection are also compared with each other.
We present two simulation designs.
In the first one, the true hazard rate is generated from a smooth age-cohort model which includes an interaction term on a small region of the age-cohort plane.
In the second case, the true hazard rate is a piecewise constant function with four heterogeneous areas.
The two true hazards are displayed in Figure \ref{fig:simu_design}, both in greyscale and in perspective plot.

The simulation design is as follows.
We set $J=10$ equally spaced age intervals $([0,10)$, $\ldots$, $[90,100])$ and $K=10$ equally spaced cohort intervals $([1900,1910)$ ,$\ldots$, $[1990,2000])$.
In order to simulate a dataset, the cohorts are first sampled uniformly over the $10$ cohort intervals and the age is then simulated using the corresponding hazard.
Censoring is then simulated as a uniform distribution over the age interval $[75, 100]$ for all cohorts such that all observed events are comprised in the age interval $[0, 100]$.
Since in practice one does not know the appropriate discretization in advance, a different discretization was used for the estimation procedure: the age and cohort intervals were defined as $5$-year length intervals instead of $10$ for the true hazard.
As a result, a total of $20\times 20$ parameters need to be estimated.
We simulated data of sample sizes $100$, $400$, $1000$, $4000$, and $10000$.
For each sample size, the simulation and estimation were replicated $500$ times.

\paragraph{Smooth true hazard}
\label{section:simu_smooth}

The smooth true hazard (Figures \ref{fig:simu_design}a and \ref{fig:simu_design}b) is generated using the age-cohort model $\log \lambda_{j,k} = \mu + \alpha_{j} + \beta_{k}$ with an intercept $\mu = \log (10 ^ {-2})$.
\textcolor{black}{The age effect vector $\bm \alpha$ and cohort effect vector $\bm \beta$ are arithmetic sequences such that $\alpha_1 = 0$, $\alpha_J = 2.5$, $\beta_1 = 0$, and $\beta_K = 0.3$.}
An interaction term is added to the hazard.
It corresponds to a bump in the hazard located in the neighbourhood of the region of the age-cohort plane (45,1945).
The bump is defined as $10$ times the Gaussian density function with mean $(1945, 45)$ and with a diagonal variance-covariance matrix with diagonal equal to $(50, 50)$.
This true hazard displays a sharp increase for high values of the age, which implies that few events will be recorded in this region.
On average, $91$~\% of the events are observed in this simulation design.
\paragraph{Piecewise constant true hazard}%
\label{section:simu_pc}

The piecewise constant true hazard (Figures \ref{fig:simu_design}c and \ref{fig:simu_design}d) has four constant areas over the age-cohort square $[0,100] \times [1900, 2000]$.
On average, $71$~\% of the events are observed in this simulation design.

\subsection{Performance of the Estimation Methods in Terms of MSE}

\begin{table}
\centering
\subfloat[Smooth true hazard]{\label{tab:mse_smooth}\scalebox{0.95}{\input{./mse_smooth.tex}}}\quad
\subfloat[Piecewise constant true hazard]{\label{tab:mse_pc}\scalebox{0.95}{\input{./mse_pc.tex}}}
\caption{Relative mean squared errors of the L$_0$ and L$_2$ methods with respect to the maximum likelihood estimate (MLE), for different sample sizes and different estimation methods. For easier comparison, the mean squared errors are given as the ratio with respect to the mean squared error of the MLE. Panel (a): smooth true hazard. Panel (b): piecewise constant true hazard.
}
\label{tab:mse}
\end{table}
Our two estimation methods (L$_{0}$ penalty and L$_{2}$ penalty) are compared in terms of the mean squared error (MSE) in each simulation scenario.
The different selection methods for the penalty (AIC, BIC, EBIC and cross-validation) are included.
\textcolor{black}{We compare our methods with the maximum likelihood estimate (MLE), which serves as baseline for comparison.
The results are presented in Table \ref{tab:mse}, which reports the relative mean square errors with respect to the MLE for easier comparison.}

Overall, the EBIC and cross-validated criteria outperform the AIC and the BIC for the two simulations scenarios.
This is particularly true for small sizes where the AIC and the BIC behave very poorly.
As expected, the L$_{2}$ penalty estimator is the most performant of all estimators in the smooth true hazard scenario (Table \ref{tab:mse_smooth}) and the L$_{0}$ method performs better in the piecewise constant hazard scenario (Table \ref{tab:mse_pc}) than in the smooth true hazard scenario.
The L$_{2}$ norm estimator is also the most performant of all estimators in the piecewise constant hazard scenario except for very large sample sizes ($n = 10000$) where the BIC, EBIC and cross-validated criterion provide slightly better performances.
In both scenarios, the EBIC always outperforms the AIC, the BIC and the cross-validated criterion.

\textcolor{black}{Finally, note that both the L$_2$ penalty and the L$_0$ penalty (with the EBIC or cross-validation) vastly outperform the MLE.
This validates that our penalized approach is useful to reduce the mean square error of the estimate. The degree to which the penalized methods outperform the MLE decreases as the sample size increases, but they still outperform the MLE for a sample size of $10000$.}

Different censoring rates were also studied which showed a degradation of the performances of the overall estimators as the percentage of censored events increases.
The performance in terms of number of selected areas was also investigated.
It showed that the EBIC and CV criterion perform better at selecting sparse models with few areas, while the AIC and BIC tend to overestimate the true number of areas.
Indeed, for sample size $4000$, the $80\%$ inter-quantile range of the selected number of areas is $\left[ 3, 5 \right]$ for the EBIC and $\left[ 1, 5 \right]$ for the CV, whereas it is $\left[ 3, 13 \right]$ and $\left[ 36, 72 \right]$ for the BIC and AIC respectively.
These experiments are not reported here.

In conclusion, the simulation experiments suggest to use the EBIC among all different criteria for the L$_{0}$ penalty as it provides the best tradeoff between computation time and estimation performance. \textcolor{black}{It has been shown that using the L$_0$ penalty is beneficial even when the true hazard is not piecewise constant, as our simulations show that the performance of this estimate exceeds that of the MLE with a smooth true hazard.}

\subsection{Perspective Plots of the Estimation Methods}

\begin{figure}
\vspace{-0.4 cm}
\begin{center}
$
\begin{array}{cc}
\includegraphics[width = 65mm, height = 65mm]{smooth_true_haz.pdf} & 
\includegraphics[width = 65mm, height = 65mm]{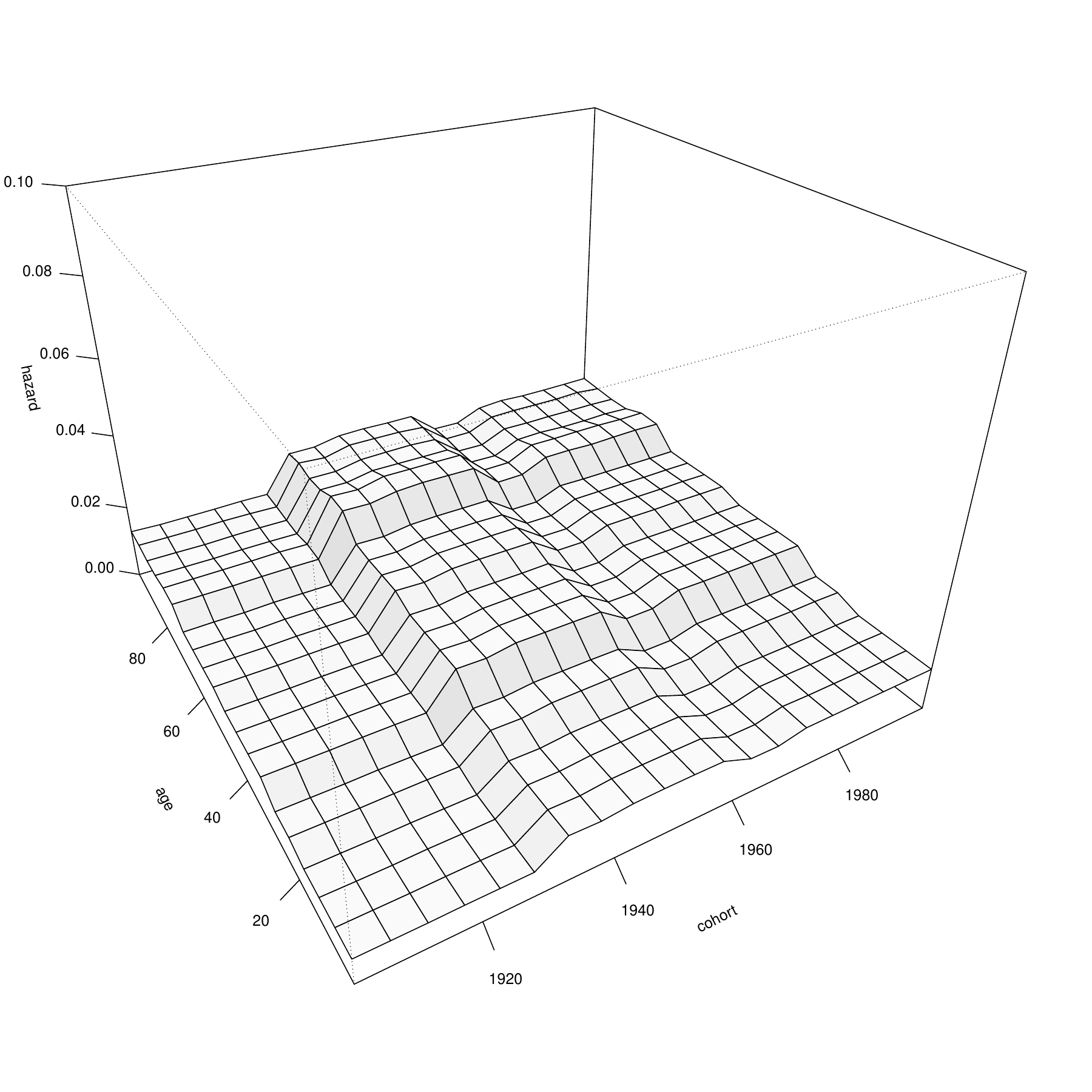} \\
\text{(a) True hazard} & \text{(b) Median of age-cohort estimates} \\
\includegraphics[width = 65mm, height = 65mm]{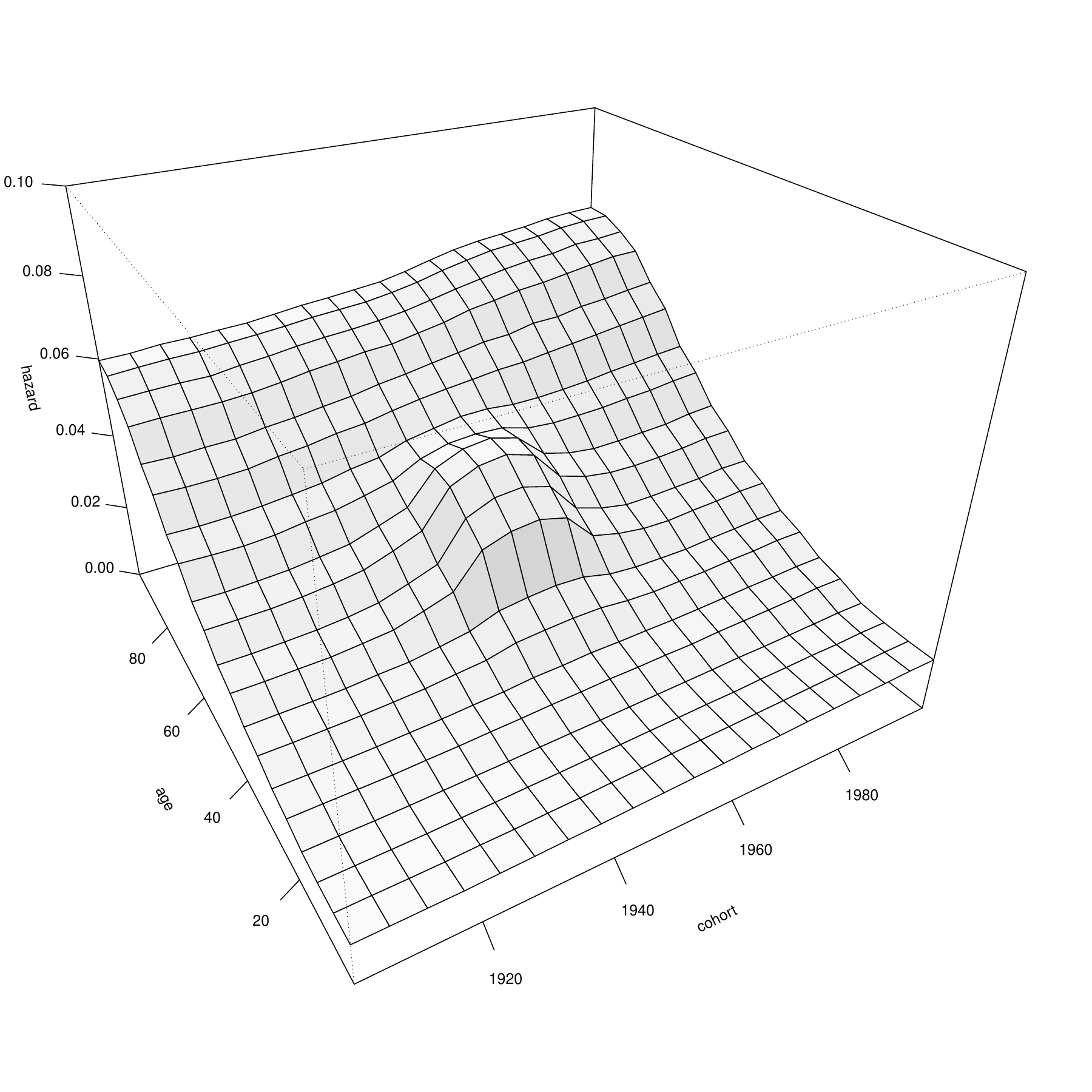} &
\includegraphics[width = 65mm, height = 65mm]{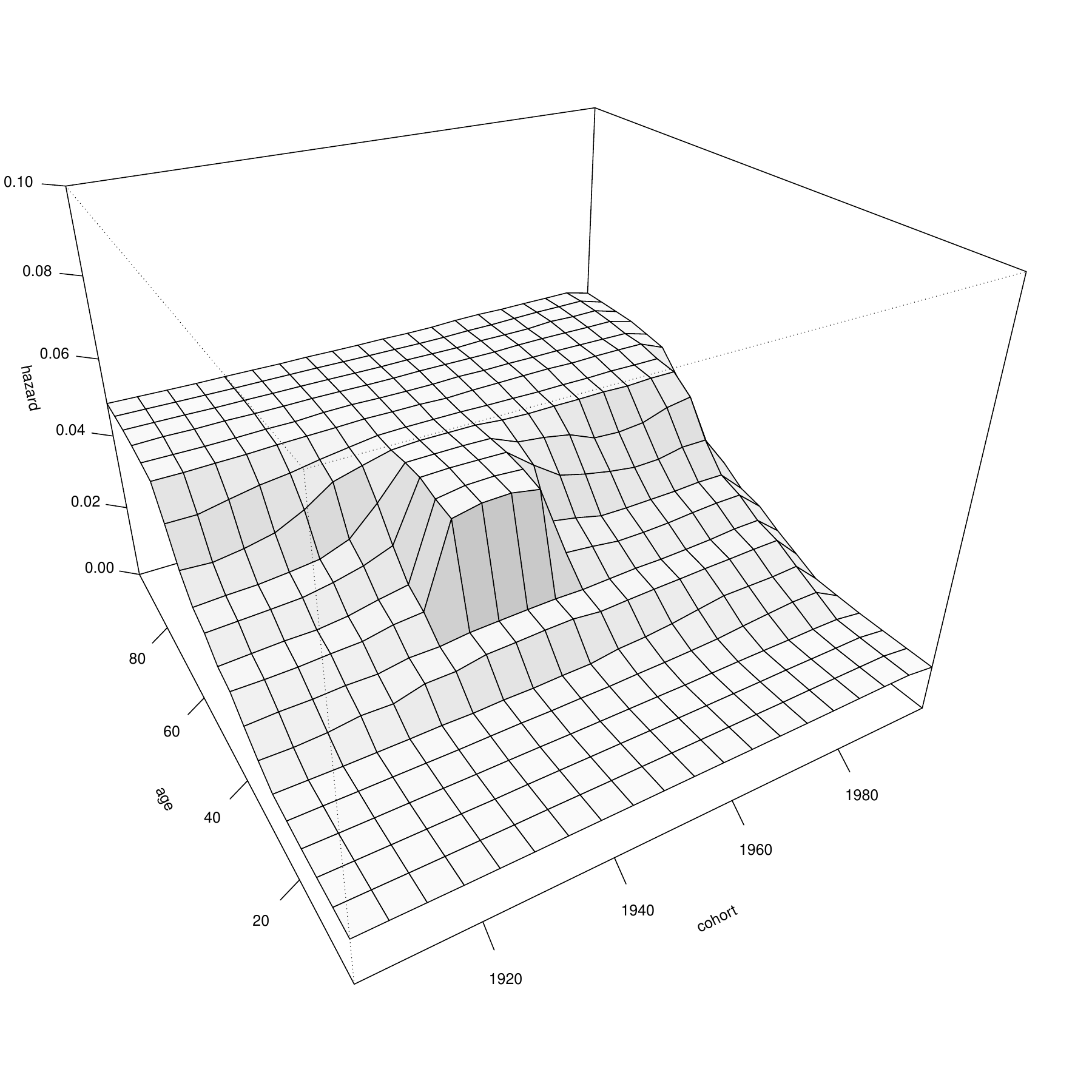} \\
\text{(c) Median of smooth estimates} & \text{(d) Median of segmented estimates}\\
\end{array}$
\end{center}
\caption{Smooth true hazard and corresponding estimates. The sample size is $4000$ and the hazard estimates are medians taken over $500$ simulations. The estimations are performed in the age-cohort plane and with different methods. Panel (a) represents the true hazard used to generate the data, Panel (b) represents the hazard estimated using the age-cohort model, Panel (c) represents the smoothed estimate, and Panel (d) represents the segmented estimate with the EBIC criterion.}
\label{fig:smooth_persp}
\end{figure}

\begin{figure}
\begin{center}
$
\begin{array}{cc}
\includegraphics[width = 70mm, height = 70mm]{pc_true_haz.pdf} & 
\includegraphics[width = 70mm, height = 70mm]{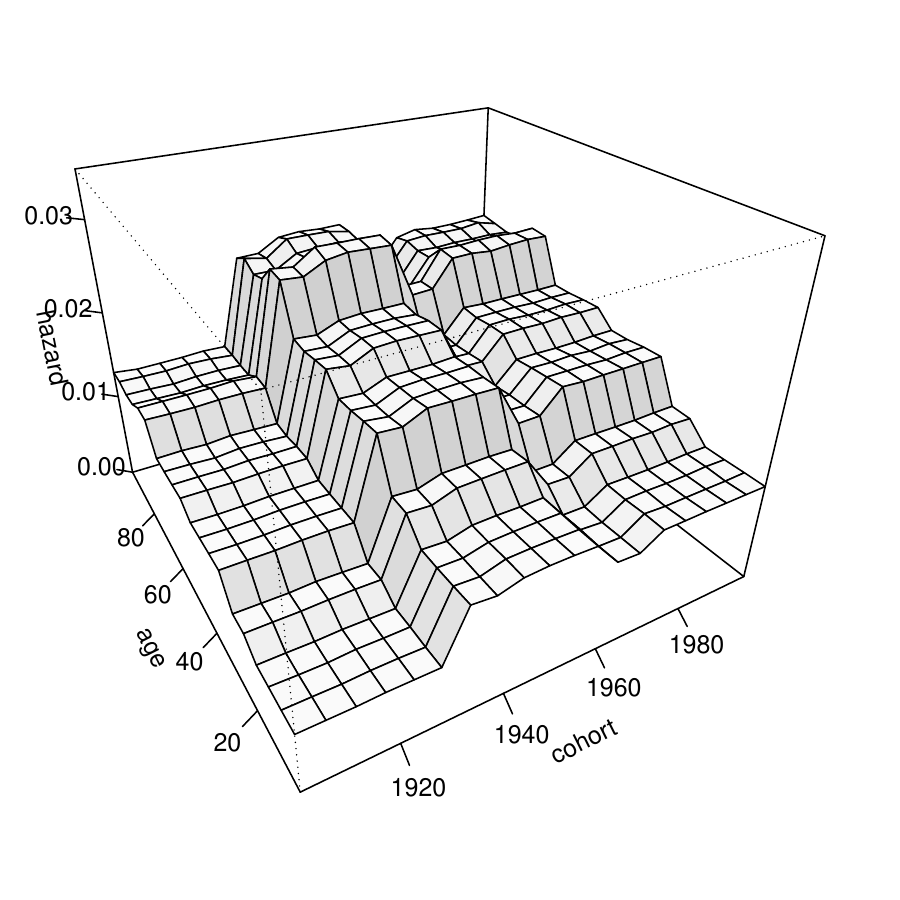} \\
\text{(a) True hazard} & \text{(b) Median of age-cohort estimates} \\
\includegraphics[width = 70mm, height = 70mm]{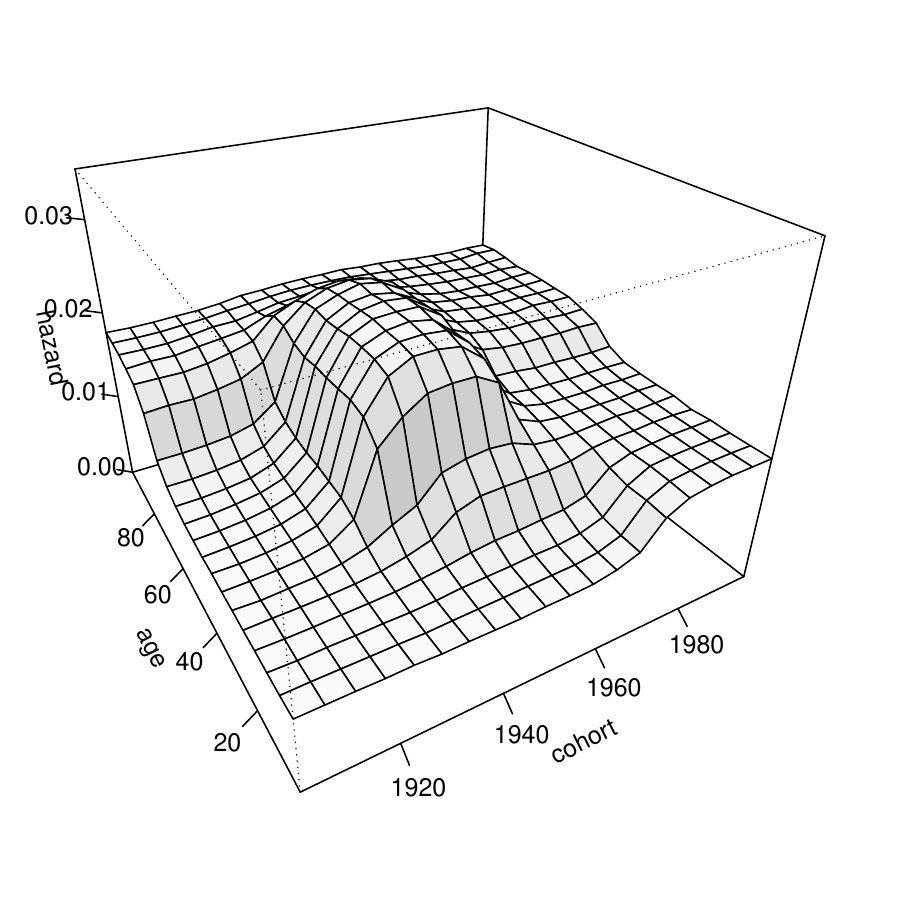} &
\includegraphics[width = 70mm, height = 70mm]{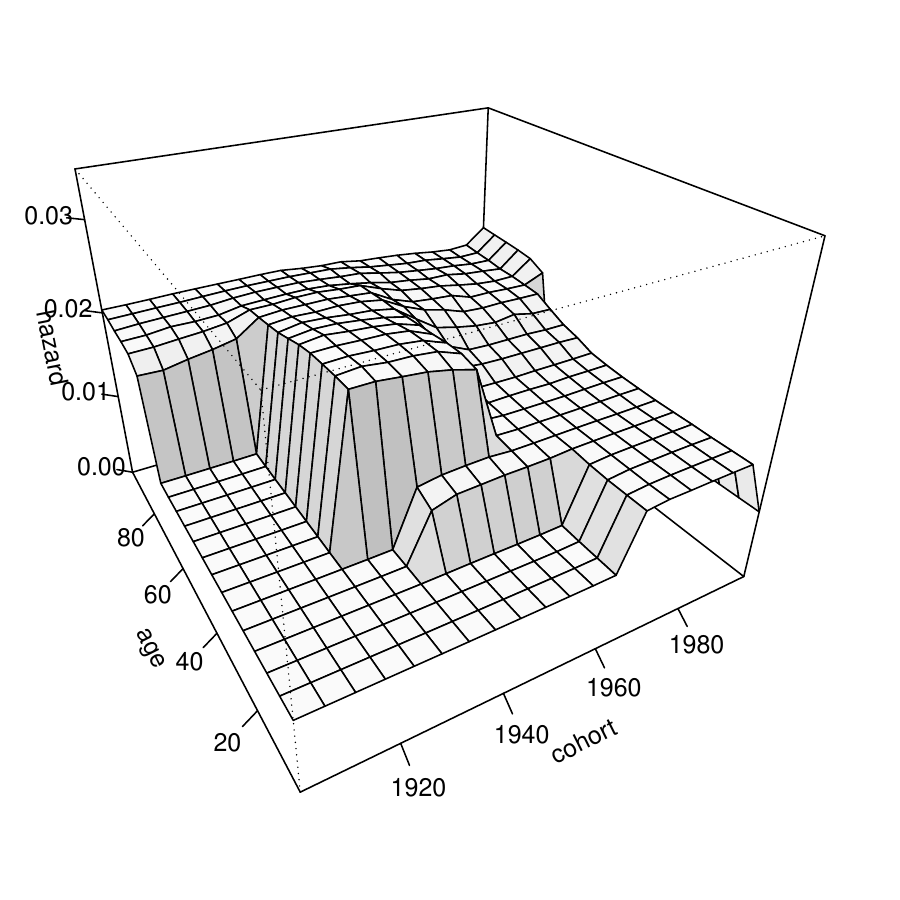} \\
\text{(c) Median of smooth estimates} & \text{(d) Median of segmented estimates}\\
\end{array}
$
\end{center}
\caption{Piecewise constant true hazard and corresponding estimates. The sample size is $4000$ and the hazard estimates are medians taken over $500$ simulations. The estimations are performed in the age-cohort plane and with different methods. Panel (a) represents the true hazard used to generate the data, Panel (b) represents the hazard estimated using the age-cohort model, Panel (c) represents the smoothed estimate, and Panel (d) represents the segmented estimate with the EBIC criterion.}
\label{fig:pc_persp}
\end{figure}

In this section the performance of our two estimates is assessed visually by comparison with the true hazard.
The standard age-cohort model \citep{Holford1983EstimationAgePeriod} has also been implemented.
This model assumes that the hazard has the following expression:
\begin{equation*}
\log \lambda_{j,k} = \mu + \alpha_{j} + \beta_{k},
\end{equation*}
where $\mu$ is the intercept, $\bm\alpha$ is the age effect and $\bm{\beta}$ is the cohort effect.
It should be noted that this model does not allow for interactions between age and cohort effects. 
Perspective plots of the median hazard estimations over $500$ replications are presented in Figures \ref{fig:smooth_persp} and \ref{fig:pc_persp} for the smooth and piecewise constant true hazard respectively.
For the L$_{0}$ regularized estimate, the penalty constant is chosen using the EBIC.

In Figure \ref{fig:smooth_persp}, it is seen that the age-cohort model is not able to estimate the central bump in the hazard. 
On the contrary, the smoothed estimate accurately recovers the shape of the true hazard except for the high values of age where few events are observed.
Interestingly, one sees that our segmentation method provides results similar to the smoothing technique even though the true hazard is not piecewise constant. 

The results in Figure 5 yield similar conclusions.
The age-cohort model behaves very poorly due to its constrained structure while the ridge and adaptive estimates provide satisfactory results.
In particular the overall shape of the true hazard is correctly estimated by the L$_0$ penalty.



\section{Real Data Application}%
\label{section:realdata}

Our method is applied to data of survival times after diagnosis of breast cancer.
The dataset is provided by the Surveillance, Epidemiology, and End Results (SEER) Program from the US National Cancer Institute (NCI).
SEER collects medical data of cancers (including stage of cancer at diagnosis and the type of tumor) and follow-up data of patients in the form of a registry.
Around $28$ percent of the US population is covered by the program.
The registry started in February $1973$ and the available current dataset includes follow-up data until January $2015$.
We refer to the website \texttt{https://seer.cancer.gov/} for information about the SEER Program and its publicly available cancer data.

In this study the duration of interest $T$ is the time from breast cancer diagnosis to death in years, the variable $U$ is the date of diagnosis (in years) and the period is the calendar time (in years).
Patients continuously entered the study between 1973 and 2015 and right-censoring occurred for patients that were still alive at the end of follow-up or for those that were lost to follow-up.

The breast cancer data was extracted using the package \texttt{SEERaBomb}.
For the sake of comparison, the subsample of malignant, non-bilateral breast tumor cancers was extracted from the dataset, such that the data comprises $1,265,277$ women with $60$ percent of censored individuals. 
Times from diagnosis to last day of follow-up vary between $0$ and $41$ years, and the dates of cancer diagnosis $U_{i}$ vary between $1973$ and $2015$.
Death from another cause than cancer is available in the dataset and is accounted for as right-censoring.

\begin{figure}
\vspace{-0.4 cm}
\begin{center}
$
\begin{array}{cc}
\includegraphics[width = 75mm, height = 75mm]{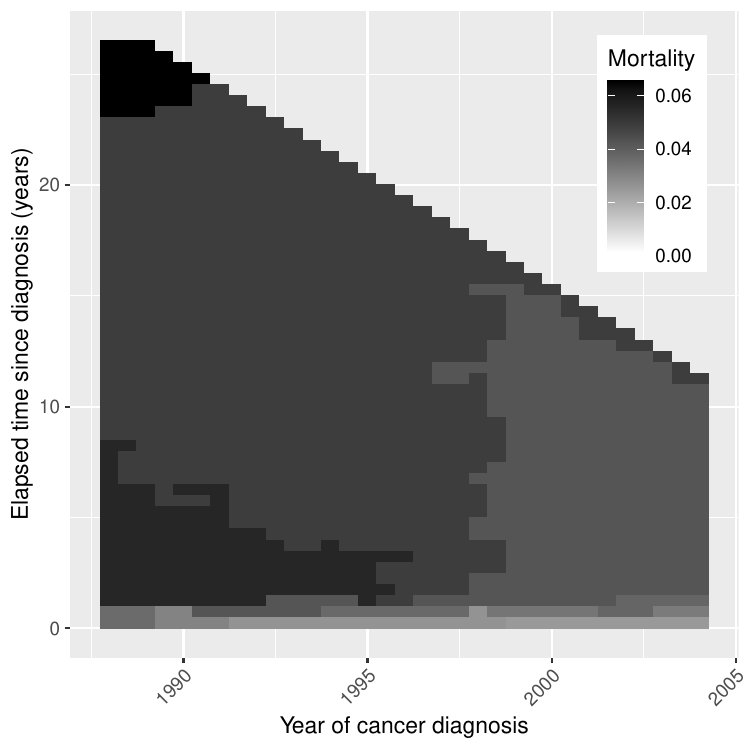} & 
\includegraphics[width = 75mm, height = 75mm]{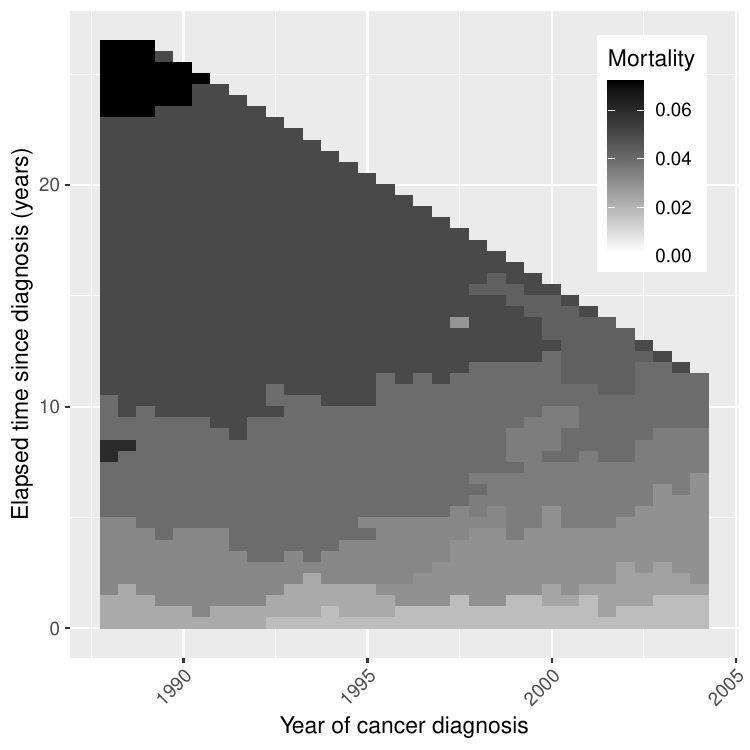} \\
\text{(a) All stages of cancer} & \text{(b) Stage 1} \\
\includegraphics[width = 75mm, height = 75mm]{stage2.pdf} & 
\includegraphics[width = 75mm, height = 75mm]{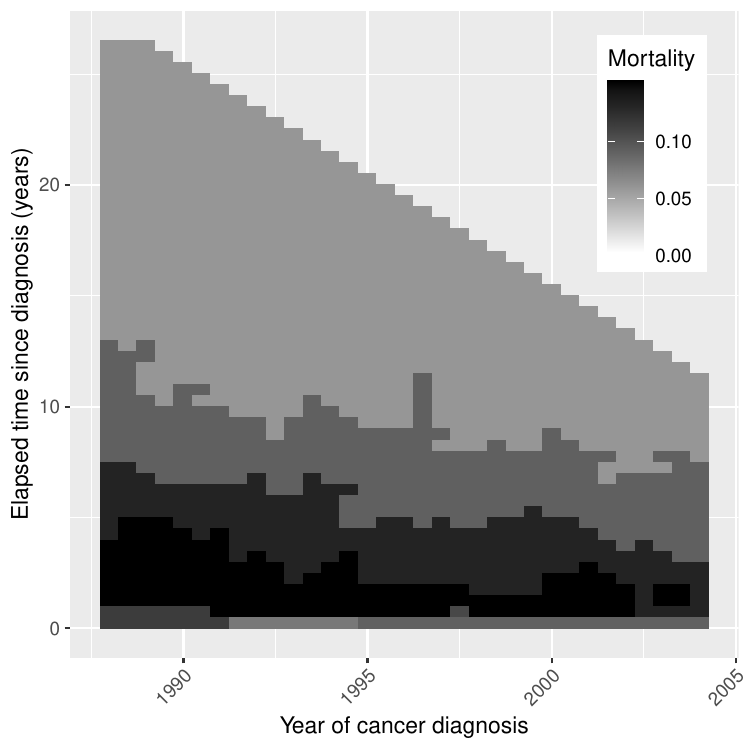} \\
\text{(c) Stage 2} & \text{(d) Stage 3} \\
\end{array}
$
\end{center}
\caption{Estimated hazard of death after diagnosis of breast cancer for different stages of cancer. The estimate is obtained with the L$_{0}$ regularization. The upper right corner of every graph corresponds to the region where no data are available. Note that the grey-color scales are different between panels.}
\label{fig:seer_stage}
\end{figure}

The implementation of our adaptive ridge method aims at two goals.
Firstly we aim at simultaneously detecting a cohort effect and an age effect, that is the evolution of the mortality with respect to the time elapsed since cancer diagnosis (age effect) and with respect to the date of diagnosis (cohort effect).
Secondly, our method will provide estimation of the hazard rates on the resulting heterogeneous areas.
The method is first applied on the whole sample of $1265277$ individuals.
In order to take into account the fact that mortality from cancer highly depends on the cancer stage, we also perform a stratified analysis with respect to the stage of cancer at diagnosis. 
For this purpose, we use the cancer stage classification provided by the SEER data: we keep the patients with cancer stages $1$, $2$, and $3$ at the time of diagnosis.
This classification closely follows that of the American Joint Comitee on Cancer (AJCC), $3^{\text{rd}}$ Edition; the details are given at page 86 of the manual entitled Comparative Staging Guide for Cancer, available at \url{https://seer.cancer.gov}.
The main difference between the two classifications is that the SEER Program classifies the cases where lymph node status cannot be assessed as if there was no regional lymph node metastasis.

The L$_{0}$ estimates for the whole sample and for each cancer stage are displayed in Figure \ref{fig:seer_stage}.
We see that the different stages of cancer at diagnosis have a great impact on the survival times.
For Stage 1 cancers, the mortality is low between $0$ and $4$--$5$ years after diagnosis, and steadily increases afterwards.
The date of diagnosis seems to have no impact on the mortality of Stage 1 cancers.
On the other hand, Stage 2 cancers exhibit a strong effect of the date of diagnosis: around $1995-1997$, the mortality significantly decreases.
This can correspond to an improvement of the treatment of breast cancer around that period in the United States.
Finally, Stage 3 cancers display a very high hazard rate across all dates of diagnosis.
This seems to indicate that the evolution in treatments of breast cancer had a significant impact on the survival times after diagnosis, but almost exclusively when cancers were diagnosed at Stage 2.
Two additional analyses of the hazard rate with stratification with respect to age at diagnosis and estrogen receptor status were performed in the Supplementary Material.
The results suggest that the shift in mortality around year $1996$ could correspond to the introduction of hormone-blocking therapy.

\section*{Conclusion}
In this article, we have introduced a new estimation method to deal with age-period-cohort analysis.
This model assumes no specific structure of the effects of age and cohort and the hazard rate is directly estimated without estimating the effects.
In order to take into account possible overfitting issues, a penalty is used on the likelihood to enforce similar consecutive values of the hazard to be equal.
Two different types of penalty terms were introduced.
One leads to a ridge type regularization while the other leads to a L$_0$ regularization.
Different selection methods of the penalty parameter were also introduced.
To our knowledge, a segmented estimation model of this kind has never been introduced in this context. 

Using simulated data, it has been shown that the cross validated ridge estimator and the $\text{EBIC}_0$ adaptive ridge estimator perform the best in terms of mean squared error.
The cross validation criterion was shown to provide the best fit of the hazard rate, but its very high computationally cost makes it non-competitive.
In this context, this modified BIC criterion comes out as a powerful tool to select the \emph{best} bias-variance tradeoff.

The method was successfully applied to data of survival after breast cancer provided by the SEER program.
The segmented estimate of the hazard rate displays important information about the shift in mortality after being diagnosed of breast cancer in the United States in the mid-$1990$s.

Our method could be directly extended to a different discretization of the age-period-cohort plane, such as $1 \times 1 \times 1$-year triangles that are represented in dark gray in Figure \ref{fig:diagram} \citep[see Section $3$ of][for an example of this discretization]{Carstensen2007AgePeriodCohort}.
Another extension would be to consider other types of penalizations.
Instead of estimating a piecewise constant hazard, one could estimate a piecewise linear hazard by penalizing over second order differences of the hazard.

\vspace*{1pc}
\noindent{\bf{Acknowledgement}}
\noindent{The authors are thankful to the National Cancer Institute for providing U.S. mortality data on cancer.}

\vspace*{1pc}
\noindent{\bf{Conflict of Interest}}
\noindent{The authors have declared no conflict of interest.}

\bibliography{biblio.bib}
\end{document}


\normalem 
\title{Supplementary Material to: Regularized Bidimensional Estimation of the Hazard Rate}

\author[1]{Vivien Goepp}
\author[1]{Jean-Christophe Thalabard}
\author[2]{Gr\'egory Nuel}
\author[1]{Olivier Bouaziz}
\affil[1]{MAP5 {\it(Department of mathematics, 45, rue des Saints-P\`eres, 75006 Paris)}}
\affil[2]{LPMA {\it{(Department of mathematics, 4, Place Jussieu, 75005 Paris)}}}
\date{January 2019}
\setcounter{Maxaffil}{0}
\renewcommand\Affilfont{\itshape\small}

\maketitle
{\color{black}
\section{Relation between the adaptive ridge and other L$_1$ and L$_2$ reweighted methods}
\label{section:mm_opt_ar}

As pointed out by a reviewer, \cite{Frommlet2016AdaptiveRidgeProcedure} provide no formal proof that the adaptive ridge approximates the L$_q$ penalty for $q\in[0, 1)$. Other iterative methods, like \cite{Candes2008EnhancingSparsityReweighted}'s L$_1$ reweighted scheme, have been shown to be an approximation of the L$_q$ penalty using a Majorization-Minimization \citep[MM, see][Section 6]{Lange2004Optimization} optimization scheme. This estimating procedure extends to the case $q=0$, where it approximates the logarithmic penalty defined below.
In this section, we show that the adaptive ridge minimizes the same function as the L$_1$ reweighted scheme. Both belong to the class of MM optimization algorithms and as such both are guaranteed to converge to a local minimum of the function to minimize.

We first prove that both methods solve the L$_q$ penalty for $0<q<1$. We then show that they both extend to the case $q=0$, where they are now approximations of the logarithm penalty instead.

\subsection{MM Optimization for Solving L$_q$ Penalties, $0<q<1$}

Consider the problem of minimizing the likelihood penalized by the L$_q$ norm:
\begin{equation} 
    \arg\min _{\bm{\beta}} \Big\{ \ell(\bm{\beta}) + \frac{\kappa}{q} \norm{\bm{\beta}}_q^q \Big\}=
    \arg\min _{\bm{\beta}} \Big\{ \ell(\bm{\beta}) + \frac{\kappa}{q} \sum_{j = 1}^p \abs{\beta_j}^q\Big\},
     \label{eq:suppl_pen_nllh}
\end{equation}
where $0<q<1$ and $\kappa>0$ is the penalty constant, rescaled here by a factor $1/q$, and $\ell(\bm{\beta})$ is the function to minimize (in our case, the negative log-likelihood).
This problem is difficult to solve because of the non-convexity of the L$_q$ norm.

We will use MM optimization to derive a numerical scheme to solving Equation \eqref{eq:suppl_pen_nllh}. MM Optimization makes use of a secondary function which majorizes the function to minimize \citep[see][]{Hunter2005Variableselectionusing}.
Since the majorization relation between functions is closed under sum, it suffices to focus in Equation \eqref{eq:suppl_pen_nllh} on the function $p(\abs{\beta_j}) =  \abs{\beta_j}_q^q /q$ for $1 \leq j \leq p$ in order to construct an MM optimization scheme.
We present two local approximations of $p(\abs*{\beta_j})$ present in the literature, which give rise to two optimization schemes.

\paragraph{L$_2$ reweighted scheme}
Let $\beta_j^{(l)} \in \mathbb{R}$ be the current point of the numerical scheme. Using a local quadratic approximation \citep[LQA, see][]{Fan2001VariableSelectionNonconcave, Hunter2005Variableselectionusing}, the function
\begin{equation*}
    q_{\text{LQA}}(\beta_j | \beta_j^{(l)}) = \frac{1}{2} \abs{\beta_j^{(l)}}^{q - 2} \beta_j ^2 + \frac{2 - q}{2q} \abs{\beta_j^{(l)}}^q
\end{equation*}
majorizes $p(\abs*{\beta_j})$ since we have $p(\abs*{\beta_j}) \leq q_{\text{LQA}}(\beta_j | \beta_j^{(l)})$ for every $\beta_j$ with equality if and only if $\beta_j = \beta_j^{(l)}$.
Define the current weights $w_j^{(l)} = \abs*{\beta_j^{(l)}}^{q-2}$.
Noting that the second term of $q_{\text{LQA}}$ does not depend on $\beta_j$, the MM optimization is given by the reweighted L$_2$ scheme:
\begin{equation}\label{eq:suppl_l2_reweighted}
\begin{split}
    \bm{\beta}^{(l)} &\gets \arg\min_{\bm{\beta}} \Big\{ \ell(\bm{\beta}) + \frac{\kappa}{2} \sum_{j = 1}^d w_j^{(l - 1)} \beta_j^2 \Big\}\\
    w_j^{(l)}&\gets \abs*{\beta_j^{(l)}}^{q - 2},
\end{split}
\end{equation}
where $(l)$ is the iteration index.

This scheme is the adaptive ridge procedure, where a small $\varepsilon$ term is added to the reweighting step to bound the denominator away from zero (see discussion on this topic hereafter).

\paragraph{L$_1$ reweighted scheme}
Using a local linear approximation \citep[LLA, see][]{Zou2008Onestepsparseestimates} the function,
\begin{equation*}
    q_{\text{LLA}}(\beta_j|\beta_j{(l)}) = \abs*{\beta_j} \abs*{\beta_j^{(l)}}^{q - 1} + \frac{1 - q}{q} \abs*{\beta_j^{(l)}}^q
\end{equation*}
majorizes $p(\abs*{\beta_j})$.
Defining now $w_j^{(l)} = \abs*{\beta_j^{(l)}}^{q - 1}$, we obtain the following reweighted L$_1$ scheme:
\begin{equation}\label{eq:suppl_l1_reweighted}
\begin{split}
   \bm{\beta}^{(l)} &\gets \arg\min_{\bm{\beta}}\Big\{ \ell(\bm{\beta}) + \kappa \sum_{j = 1}^d w_j^{(l - 1)} \abs*{\beta_j} \Big\}\\
   w_j^{(l)}&\gets \abs*{\beta_j^{(l)}}^{q - 1}.
\end{split}
\end{equation}

\subsection{Extension to the case $q = 0$}

Let us note that even though $q$ has to be strictly positive in Equation \eqref{eq:suppl_pen_nllh}, the numerical schemes \eqref{eq:suppl_l2_reweighted} and \eqref{eq:suppl_l1_reweighted} are still defined for $q = 0$.
We now show that in the case $q=0$, both schemes do not solve the L$_0$ penalty: they correspond to a logarithmic penalty, which is a good approximation thereof \citep[][Section 2.3]{Candes2008EnhancingSparsityReweighted}.
\textcolor{black}{Let us first note that formally, the logarithmic penalty seems to be a good approximation to the L$_0$ norm since $\lim_{q \to 0} \sum_{j = 1}^p \abs*{\beta_j}^q \to \norm*{\bm\beta}_0$ and $\lim_{q\to0} (1 / q) \sum_{j = 1}^p (\abs*{\beta_j}^q - 1) = \sum_{j = 1}^p\log(\abs*{\beta_j})$.
In the context of sparse signal recovery, this is enough to prove that the logarithmic penalty yields the same estimate as the L$_0$ penalty (including the case with $\varepsilon$ perturbation) \citep[Section I]{Wipf2010Iterative}, although this does not seem to be easy to prove in the case of penalized likelihood.} 

We will start with the case of the L$_2$ reweighted scheme.
Consider Problem \eqref{eq:suppl_pen_nllh} where the L$_q$ penalty is replaced by the logarithmic penalty: $p(\abs*{\beta_j}) = \log(\abs*{\beta_j})$.
The LQA of this penalty around the current point $\beta_j^{(l)}$ is given by
\begin{align*}
    q(\beta_j|\beta_j^{(l)})& = p(\abs*{\beta_j^{(l)}}) + \left(\beta_j^2 - \big(\beta_j^{(l)}\big)^2\right) \frac{p'(\abs*{\beta_j^{(l)}})}{2\abs*{\beta_j^{(l)}}}\\
    & = \frac{1}{2} \frac{\beta_j^2}{\big(\beta_j^{(l)}\big)^2} + \log(\abs*{\beta_j^{(l)}}) - \frac{1}{2}
\end{align*}
and the MM optimization scheme is obtained by iteratively minimizing $q(\beta_j|\beta_j^{(l)})$:
\begin{equation}
\begin{split}\label{eq:letter_ar_mm}
\bm{\beta}^{(l)} & \gets \arg\min_{\bm{\beta}} \Big\{\ell_n(\bm{\beta}) + \frac{\kappa}{2} \sum_{j=1}^p \frac{\beta_{j}^2}{\big(\beta_j^{(l)}\big)^2} \Big\}, \\
\end{split}
\end{equation}
which is the adaptive ridge with $q = 0$, with $\varepsilon$ set to zero.
It is straightforward to show that the case $\varepsilon>0$ corresponds to the penalty function $\beta_j \mapsto \log(1 + \beta_j^2 / \varepsilon^2)$ instead, which is defined for $\beta_j = 0$. 
The theoretical properties of the $\varepsilon$-perturbed LQA is studied in \cite{Hunter2005Variableselectionusing} for a specific class of penalties.

The same reasoning applies to the LLA and shows \citep{Zou2008Onestepsparseestimates, Candes2008EnhancingSparsityReweighted} that the L$_1$ reweighted scheme with $q = 0$ corresponds to the MM optimization of~\eqref{eq:suppl_pen_nllh} with penalty function $p(\abs*{\beta_j})= \log(\abs*{\beta_j} / \abs*{\varepsilon} + 1).$

\subsection{Related Works}

Many works have made use of the reweighted L$_1$ and L$_2$ schemes derived above. We mention some related works of importance and finish with some remarks on the relative merits of the two methods.

These methods seem to first have been used in compressed sensing: \cite{Candes2008EnhancingSparsityReweighted} studied the L$_1$ reweighted scheme with $q = 0$, while \cite{Daubechies2010Iterativelyreweightedleast} and \cite{Chartrand2008Iterativelyreweightedalgorithms} studied Algorithm \eqref{eq:suppl_l2_reweighted} for $0<q\leq1$ and $q\in[0,1)$ respectively.
\cite{Johnson2012LogPenalizedLeastSquares} studied Algorithm \eqref{eq:suppl_l1_reweighted} with $q = 0$ in the context of linear regression.
\cite{deRooi2011Deconvolutionpulsetrains}, \cite{Rippe2012VisualizationGenomicChanges}, and \cite{deRooi2014SparseDeconvolutionOne} used Algorithm \eqref{eq:suppl_l2_reweighted} with $q=0$ in various applications, while \citet[][Section 5]{Bach2011OptimizationSparsityInducingPenaltiesa} and \citet[][Section 5.4]{Mairal2014Sparse} used the LQA to derive Algorithm \eqref{eq:suppl_l2_reweighted} for the L$_q$ penalty ($0<q<1$) and for more general norms.
More recently, \cite{Frommlet2016AdaptiveRidgeProcedure} studied Algorithm \eqref{eq:suppl_l2_reweighted} numerically for $q\in[0, 1)$, and specifically $q = 0$, under the name ``adaptive ridge'', which is the method used in this work. \cite{Dai2018Brokenadaptiveridge} proved its consistency and oracle property in the setting of linear regression.
Finally, \cite{Tardivel2018Sparsestrepresentationsapproximations} have recently proven that, in the case of sparse signal recovery, the L$_1$ reweighted scheme with $q=0$ is equivalent to minimizing the L$_0$ penalty problem.

\paragraph{Remark 1.} The choice of $q$ is independent from the choice between reweighted L$_1$ or L$_2$ schemes and is not tackled here. Many papers cited in this section seem to favor choosing a small value of $q$.
\paragraph{Remark 2.} Both reweighted L$_1$ and L$_2$ schemes have their advantages and drawbacks. The former is sparse at every step but each step requires solving a L$_1$ penalty. The latter is only asymptotically sparse and thus may require more iterations but it involves the simpler L$_2$ penalty, whose solution is explicit in the linear regression setting and simple to derive in other settings. To the best of our knowledge, there is no available implementation of the fused L$_1$ penalty for a general negative log-likelihood $\ell_n(\bm{\beta})$ .
\paragraph{Remark 3.} As in the present work, most works cited in this section use a modified weighting step for numerical stability: the denominator is bounded away from zero with an $\varepsilon$ perturbation. While some offer rules of thumb to adaptatively decrease the value of $\varepsilon$ as the algorithm converges, we have followed \cite{Frommlet2016AdaptiveRidgeProcedure}'s implementation and have set $\varepsilon$ to a very small fixed value.
}

\section{Application to Breast Cancer Mortality: Stratification with Respect to the Age at Diagnosis}%
\label{section:age_diagnosis}

\begin{figure}
		\centering
		\includegraphics[width = 170 mm, height = 150 mm, keepaspectratio]{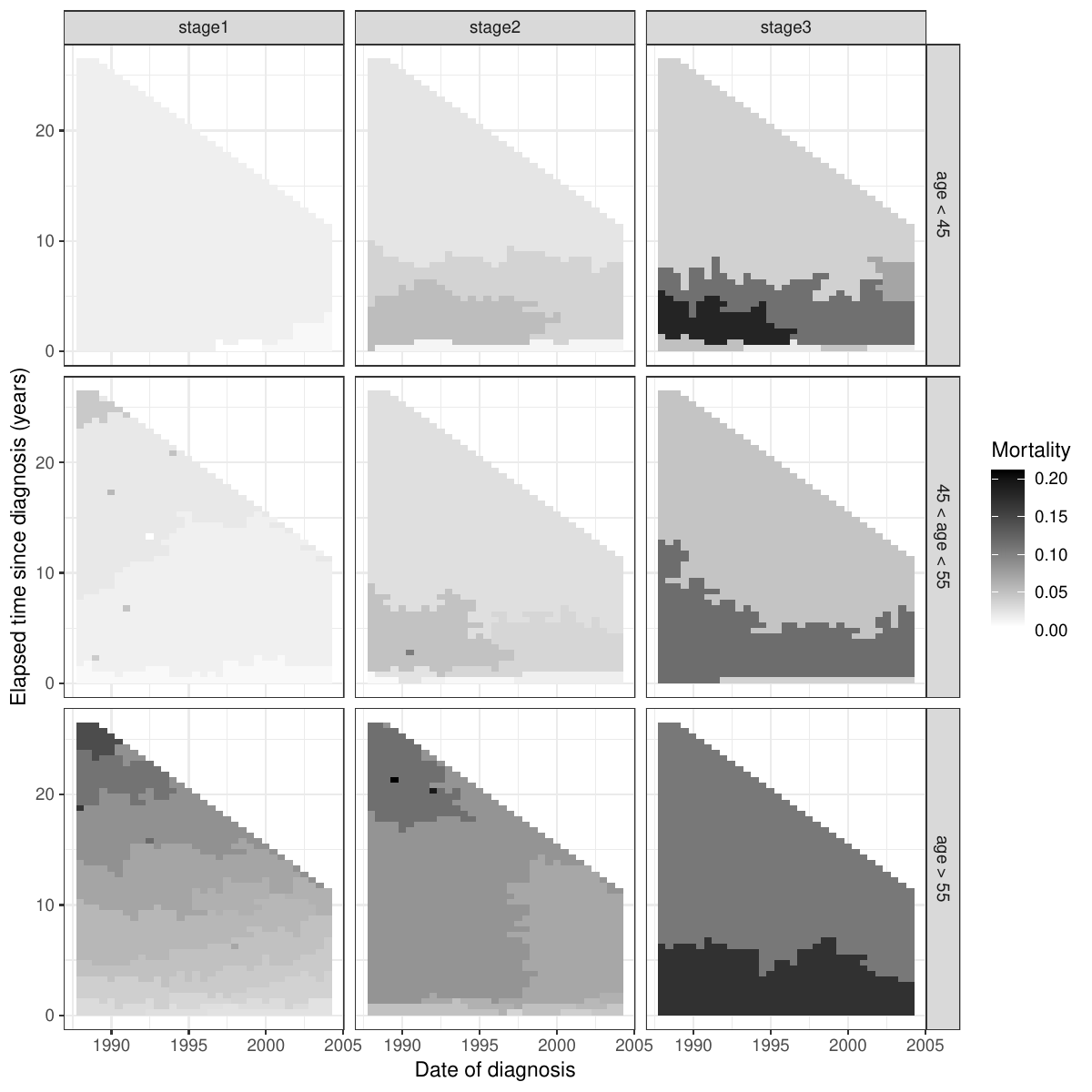}
\caption{Estimated hazard of death since diagnosis of breast cancer for different cancer stages and for different ages at diagnosis. The estimate is obtained with the L$_{0}$ regularization. The upper right corner of every graph corresponds to the region where no data are available. All graphs share the same scale.}
\label{fig:seer_stage_age}
\end{figure}

The mortality of breast cancer is known to greatly vary on whether the cancer is pre or post-menopausal \citep{Consensus1985}.
Consequently, a thorough analysis of the mortality from breast cancer would require to stratify with respect to the menopausal status at diagnosis.
Since this covariate is not present in the data, we decided to stratify the sample with respect to the age of the patient at diagnosis, which is a proxy of menopausal status.
Most women are known to have their menopause between $45$ and $55$ years old \citep{Hill1996, Henderson2008, Gold2011}, with 25th, 50th, 75th percentiles ranging from years 47-49, 50-51, 52-54, respectively, according to countries and surveys \citep{Mishra2017}.
Consequently, based on the available information in SEER, for each cancer stage, the patients were divided into three classes of age at diagnosis: $(., 45]$, $(45, 55]$, and $(55, .)$ as a proxy for pre- menopausal, peri- menopausal and post- menopausal ages, respectively. 
The resulting estimated hazards are represented in Figure~\ref{fig:seer_stage_age}.

The stage $1$ cancer patients younger than $45$ and the stage $3$ cancer patients older than $55$ display the same mortality across all dates of diagnosis, i.e. with no cohort effect.

Moreover, the mortality of stage $1$ cancer patients aged $45$ and older at diagnosis has a slight cohort effect corresponding to a progressive decrease in the mortality across all survival times \citep{Peto2000}.
This could suggest a trend of slow and steady improvement of the treatment of breast cancer in the United States over the period $1887-2005$.

Finally, we observe a clear decrease of the mortality for stage $2$ cancers for all three age classes.
This shift is located at the year $1995$ for middle-aged patients and around the years $1997-1998$ for patients younger than $45$ and older than $55$.
The same drop in mortality is observed for stage $3$ cancers with patients younger than $45$ at diagnosis, around year $1995$.
This could correspond to the introduction of improvements in the treatments of breast cancer in the United States \citep{Consensus1985}.
Among the three main medical innovations, which can be considered in this period, the improvement of the surgical procedures for the loco- regional control of the disease and the assessment of the beneficial effect of hormone-receptor therapies could be reflected in the observed survival in stages 1-2, whereas the later emergence during this period of new classs of chemotherapeutic agents like taxoids \citep{Rowinsky1992, Crown2004} or herceptin-based therapies targeted on new class of tumor markers \citep{Pegram1998, Emens2004} would be related with the changes in survival observed in stage 3. 
In the next section, we will use a stratified analysis to understand the effect of hormone-receptor therapies on the mortality shift in the mid-$1990$s.

\section{Application to Breast Cancer Mortality: Stratification with Respect to the Estrogen Receptor Status}%
\label{section:er_status}

\begin{figure}
		\centering
		\includegraphics[width = 200 mm, height = 180 mm, keepaspectratio]{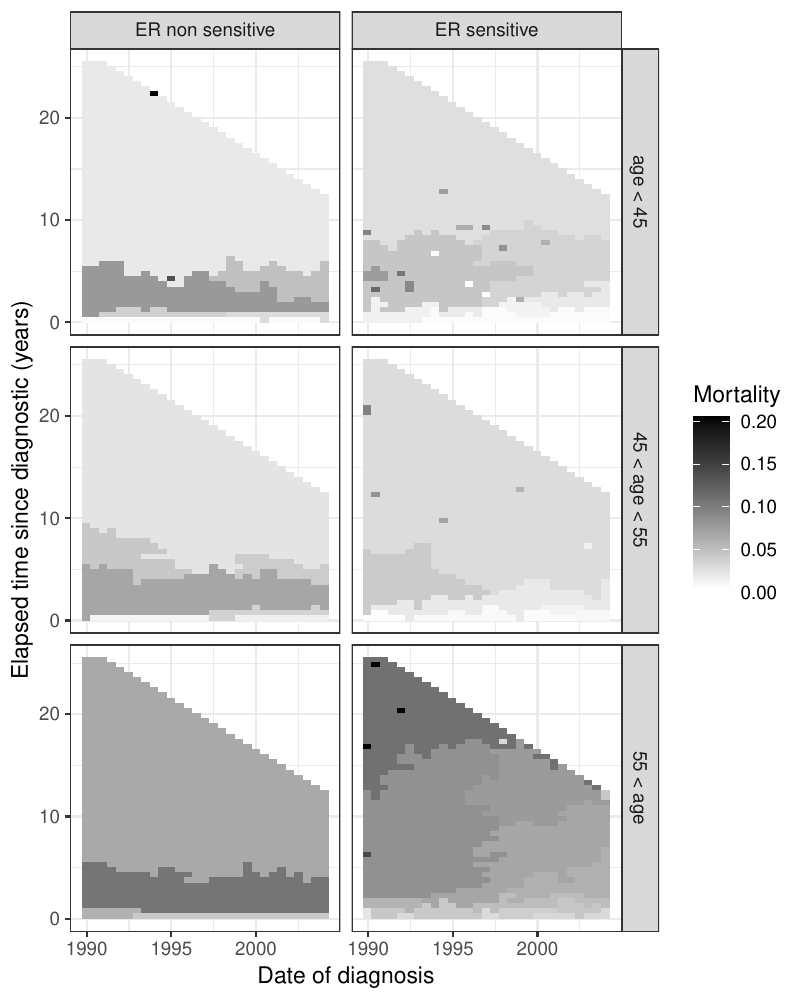}
\caption{Estimated hazard of death since diagnosis of breast cancer for Stage $2$ cancers. The estimation is carried separately for three classes of age at diagnosis: $(., 45]$, $(45, 55]$, and $(55, .)$ and for sensitive and non-sensitive estrogen receptor cancers. Inference is made with the L$_{0}$ regularization. All graphs share the same scale.}
		\label{fig:stage2_er}
\end{figure}

The cohort effect highlighted in the previous section could correspond to the introduction of Selective Estrogen Receptor Modulator (SERM) treatments and in particular the use of Tamoxifen as a treatment for breast cancer, showing improved survival in women with estrogen receptor positive tumor, initially in post- menopausal women \citep{Fisher1989}, later in both post - and pre- menopausal women \citep{Group1988, Fisher1998TamoxifenPreventionBreasta, Pritchard2005, Cochrane2008}.
Indeed, Tamoxifen was gradually used in the early years of $1990$'s \citep{Gail1999WeighingRisksBenefits, Harlan2002, Mariotto2006} to decrease the mortality of breast cancer patients.
This treatment is only efficient on estrogen receptor-sensitive cancers.
To validate our hypothesis, we conducted the estimation of mortality separately for patients with estrogen receptor sensitive and non-sensitive cancers.
Since stage $2$ cancers displayed a strong cohort effect across all ages at diagnosis, we only kept stage $2$ cancers in this study.
The estimated mortality is given in Figure~\ref{fig:stage2_er}.
Note that the spikes in the mortality are an artifact of the segmentation procedure when the sample sizes tend to be too small in some regions of the age-cohort plane and are not to be taken into account in the interpretation of the mortality.

There is a clear difference in the evolution of mortality with respect to time at diagnosis between sensitive and non-sensitive estrogen cancers.
For estrogen sensitive cases, the mortality displays the same sudden decrease around years $1997-1998$ as in Figure~\ref{fig:seer_stage_age}, across all age classes.
In particular for individuals aged $55$ or more at the time of diagnosis, the mortality has gradually decreased for estrogen sensitive patients, whereas it did not evolve with time for estrogen non-sensitive patients.
On the other hand, the mortality for non-estrogen sensitive cancers displays almost no cohort effect for all ages at diagnosis \citep{Knight1977}.

The same analysis was run with stratification with respect to progesterone receptor status, with very similar morality estimates (results not shown here).
Further analyses could be carried out to better understand the effect of the introduction of hormone-blocking therapies on mortality.
However, the segmentation of the hazard rate, even with this simple stratified analysis, highlighted that the adoption of SERM therapies in the United States is a potential reason for the sharp decrease of mortality in the middle of the $1990$s \citep{Peto2000}.

\bibliographystyle{abbrvnat}
\bibliography{biblio.bib}

%% file: mse_smooth.tex
\scalebox{0.9}{
\begin{tabular}{lrrrrrr}
  \toprule
   & \multicolumn{4}{c}{\textbf{L$_{0}$ method}} & \multicolumn{1}{c}{\textbf{L$_{2}$ method}} & \multicolumn{1}{c}{\textcolor{black}{\textbf{MLE}}}\\
 \cmidrule(lr){2-5}
 \textbf{Sample size} & \textbf{AIC} & \textbf{BIC} & \textbf{EBIC} & \textbf{CV} & \textbf{CV} & \textbf{}\\
 \midrule
100 & 1.016 & 0.988 & 0.011 & 0.011 & 0.002 & 1 \\ 
  400 & 1.005 & 0.845 & 0.144 & 0.026 & 0.004 & 1 \\ 
  1000 & 0.946 & 0.628 & 0.024 & 0.020 & 0.006 & 1 \\ 
  4000 & 0.851 & 0.267 & 0.054 & 0.037 & 0.011 & 1 \\ 
  10000 & 0.634 & 0.144 & 0.113 & 0.057 & 0.024 & 1 \\ 
   \bottomrule
\end{tabular}
}

%% file: mse_pc.tex
\scalebox{0.9}{
\begin{tabular}{lrrrrrr}
  \toprule
   & \multicolumn{4}{c}{\textbf{L$_{0}$ method}} & \multicolumn{1}{c}{\textbf{L$_{2}$ method}} & \multicolumn{1}{c}{\textcolor{black}{\textbf{MLE}}}\\
 \cmidrule(lr){2-5}
 \textbf{Sample size} & \textbf{AIC} & \textbf{BIC} & \textbf{EBIC} & \textbf{CV} & \textbf{CV} & \textbf{}\\
 \midrule
100 & 1.004 & 1.001 & 0.003 & 0.003 & 0.002 & 1 \\ 
  400 & 0.984 & 0.775 & 0.036 & 0.029 & 0.012 & 1 \\ 
  1000 & 0.829 & 0.408 & 0.092 & 0.085 & 0.024 & 1 \\ 
  4000 & 0.715 & 0.128 & 0.090 & 0.110 & 0.058 & 1 \\ 
  10000 & 0.720 & 0.083 & 0.065 & 0.081 & 0.107 & 1 \\ 
   \bottomrule
\end{tabular}
}